\documentclass[a4paper,11pt]{article}

\usepackage[utf8]{inputenc}
\usepackage{fullpage}
\usepackage{amsmath}
\usepackage{amsthm}
\usepackage{amssymb}
\usepackage{hyperref}
\usepackage{bbm}
\usepackage{url}
\usepackage{graphicx}
\usepackage{tikz-cd}
\usetikzlibrary{decorations.pathmorphing}
\usepackage{stmaryrd}
\usepackage[titletoc,title]{appendix}
\usepackage[nottoc]{tocbibind}

\usepackage{enumitem}
\makeatletter
\newcommand{\mylabel}[2]{#2\def\@currentlabel{#2}\label{#1}}
\makeatother

\usepackage[framemethod=TikZ]{mdframed}

\numberwithin{equation}{section}
\usepackage[margin=0.5cm, font=small]{caption}
\tikzcdset{every label/.append style = {font = \normalsize}}
\theoremstyle{definition}
\newtheorem{dfn}{Definition}[section]
\newtheorem{prop}[dfn]{Proposition}
\newtheorem{thm}[dfn]{Theorem}
\newtheorem{lem}[dfn]{Lemma}
\newtheorem{corol}[dfn]{Corollary}
\newtheorem{fact}[dfn]{Fact}
\newtheorem*{rem*}{Remark}

\newenvironment{fprop}
{\begin{mdframed}[roundcorner=10pt]
 \begin{prop}}
{\end{prop}
 \end{mdframed}}

\newenvironment{flem}
{\begin{mdframed}[roundcorner=10pt]
 \begin{lem}}
{\end{lem}
 \end{mdframed}}

\newenvironment{fthm}
{\begin{mdframed}[roundcorner=10pt]
 \begin{thm}}
{\end{thm}
 \end{mdframed}}

\newenvironment{fcorol}
{\begin{mdframed}[roundcorner=10pt]
 \begin{corol}}
{\end{corol}
 \end{mdframed}}

\newenvironment{ffact}
{\begin{mdframed}[roundcorner=10pt]
 \begin{fact}}
{\end{fact}
 \end{mdframed}}

\newcommand{\shadowcolor}{gray!20}

\newenvironment{shadow}
{\begin{mdframed}[linewidth=0pt, backgroundcolor=\shadowcolor]}
{\end{mdframed}}

\newenvironment{sproof}
{\begin{shadow}\begin{proof}}
{\end{proof}\end{shadow}}

\makeatletter
\newcommand{\dotminus}{\mathbin{\text{\@dotminus}}}
\newcommand{\@dotminus}{%
  \ooalign{\hidewidth\raise1ex\hbox{.}\hidewidth\cr$\m@th-$\cr}%
}
\makeatother

\newcommand{\N}{\mathbb{N}}

\newcommand{\one}{\mathbbm{1}}

\newcommand{\arrowset}[1]{\overset{#1}{\longrightarrow}}
\newcommand{\arrowsetl}[1]{\overset{#1}{\longleftarrow}}

\newcommand{\Id}{{\operatorname{Id}}}

\newcommand{\len}{{\operatorname{len}}}
\newcommand{\tr}{{\operatorname{tr}}}
\newcommand{\tail}{{\operatorname{tail}}}
\newcommand{\zerothDef}{{\operatorname{zerothDef}}}
\newcommand{\nthDef}{{\operatorname{nthDef}}}
\newcommand{\nth}{{\operatorname{nth}}}
\newcommand{\deff}{{\operatorname{def}}}
\newcommand{\IdUntil}{{\operatorname{IdUntil}}}
\newcommand{\List}{{\operatorname{List}}}
\newcommand{\Seq}{{\operatorname{Seq}}}
\newcommand{\concat}{\mathbin{+\mkern-2mu+}}

\newcommand{\Set}{\textbf{Set}}

\newcommand{\param}[1]{parametrized}

\newcommand\blfootnote[1]{%
  \begingroup
  \renewcommand\thefootnote{}\footnote{#1}%
  \addtocounter{footnote}{-1}%
  \endgroup
}

\begin{document}

\title{The List Object Endofunctor is Polynomial}
\date{\today}
\author{By Samuel Desrochers}
\maketitle

\begin{abstract}
In this paper, we study the list object functor $L : \mathcal{C} \rightarrow \mathcal{C}$ for a general category $\mathcal{C}$ with finite limits and \param{} list objects. We show that $L$ is polynomial as long as $\mathcal{C}$ is extensive.
\end{abstract}

\tableofcontents


\pagebreak
\section{Introduction}
\label{sect:intro}

\subsection{Intuition and main result}
\label{sect:intro:intuit}

Given a set $X$, it is possible to form the set $L(X)$ of finite lists with elements in $X$, also called \textit{lists on $X$}. This mapping $X \mapsto L(X)$ extends to a functor $L : \Set \rightarrow \Set$, where $f : X \rightarrow Y$ becomes the map $L(f) : L(X) \rightarrow L(Y)$ given by $L(f)([x_1, ..., x_n]) = [f(x_1), ..., f(x_n)]$.

One way of thinking of the collection $L(X)$ of lists on $X$ is to divide lists by length. For any fixed $n \in \N$, the set of lists on $X$ of length $n$ can just be thought of as $X^n$, the product of $X$ with itself $n$ times. This gives us a new way to express the set $L(X)$:
\begin{align*}
L(X) = \coprod \limits_{n \in \N} X^n.
\end{align*}
This expression, a ``sum of powers of $X$", is reminiscent of a polynomial, and indeed, $L$ is an example of a \textit{polynomial functor}. This means that $L$ is ``represented" by a \textit{polynomial}, a certain type of diagram (we recall these definitions in section \ref{sect:prelim}). In our case, the polynomial in question is
\begin{equation*}
\begin{tikzcd}
\one & E \arrow[l] \arrow[r, "\pi_2^E"] & \N \arrow[r] & \one,
\end{tikzcd}
\end{equation*}
where $E = \{(m,n) \in \N \times \N : m<n\}$, and $\pi_2^E : E \rightarrow \N$ is the second projection.

The goal of this paper is to state and prove that $L$ is a polynomial functor in a much more general setting. First, the idea of ``lists on $X$" can be extended to categories more general than $\Set$ by using \textit{list objects}, introduced by Cockett \cite{CockettList} (see also section \ref{sect:NNOlist:list}). In a category $\mathcal{C}$ with finite limits where every object $X$ has an associated (\param{}) list object $L(X)$, we can again form $L : \mathcal{C} \rightarrow \mathcal{C} : X \mapsto L(X)$ and ask whether it is a polynomial functor. To show this is the case, we need to use an analogue of the above diagram: $\N$ will be replaced by the \textit{natural numbers object} $N = L(\one)$ (see section \ref{sect:NNOlist:NNO}), and $E$ will be constructed using finite limits (see section \ref{sect:NNOlist:E}). Finally, we will need to assume that the category $\mathcal{C}$ is \textit{extensive} (we recall the definition in section \ref{sect:prelim}); in the terminology of \cite{CockettList}, this means $\mathcal{C}$ is a \textit{locos}.

In short, the goal of this paper is to prove the following theorem. The proof is found in section \ref{sect:outline}, along with an extended introduction and a proof outline; see also section \ref{sect:intro:structure} for the paper structure.

\begin{fthm}
\label{thm:conclusion}
Let $\mathcal{C}$ be an extensive category with finite limits and \param{} list objects. Then the list object functor $L : \mathcal{C} \rightarrow \mathcal{C}$ is represented by the polynomial
\begin{equation*}
\begin{tikzcd}
\one & E \arrow[l] \arrow[r, "\pi_2^E"] & N \arrow[r] & \one.
\end{tikzcd}
\end{equation*}
\end{fthm}

The functor $L$ has been shown to be part of a Cartesian monad \cite[corollary 2.7]{JayCartMonad}, so this means that $L$ is part of a \textit{polynomial monad}.

In a category with finite limits, a polynomial of the form $I \arrowsetl{s} A \arrowset{f} B \arrowset{t} J$ has an associated polynomial functor as long as $f$ is exponentiable (see section \ref{sect:prelim}), so all polynomial functors exist in a locally cartesian closed category (see \cite{Polynomials}). Theorem \ref{thm:conclusion} is notable because it does not assume that $\mathcal{C}$ has any exponentials: the fact that $\pi_2^E : E \rightarrow N$ is exponentiable follows from the existence of list objects.

\subsection{Structure of this paper}
\label{sect:intro:structure}

The first few sections of this paper are preliminaries: section \ref{sect:prelim} reviews standard categorical concepts, section \ref{sect:IL} presents the internal language that is used throughout the paper, and section \ref{sect:NNOlist} covers natural numbers objects and list objects (using the internal language), including some important constructions used in the paper, such as $E$ and $\pi_2^E$.

Then, in section \ref{sect:outline}, we explain how to prove theorem \ref{thm:conclusion}. We show how to reduce this theorem to a more tangible result, theorem \ref{thm:finalPolyn}, whose proof takes up the rest of the paper (the proof of theorem \ref{thm:conclusion} is in section \ref{sect:outline}, contingent on theorem \ref{thm:finalPolyn}). We also note that the assumption of extensivity for $\mathcal{C}$ can be replaced by two weaker hypotheses, \ref{hyp:list} and \ref{hyp:NNO}.

Finally, sections \ref{sect:technical} through \ref{sect:equality} contain the main elements for the proof of theorem \ref{thm:finalPolyn}; section \ref{sect:outline} explains how these parts fit together. Note that, throughout sections \ref{sect:technical} to \ref{sect:equality}, we assume that we're working in a (not necessarily extensive) category $\mathcal{C}$ with finite limits and \param{} list objects (this is not always restated at each result).

\subsection{Acknowledgements}

This research was conducted as part the author's PhD thesis, under the supervision of Simon Henry and Philip J. Scott. The author wishes to express his gratitude for their insights and support.

This research was partially supported by the Natural Sciences and Engineering Research Council of Canada (NSERC) and by the Ontario Graduate Scholarship (OGS) Program.

\section{Standard categorical notions}
\label{sect:prelim}

In this section, we'll review standard concepts from category theory to establish notation and to restate useful results from the literature.

\subsubsection*{Products and pullbacks}

A category which has all finite products is called \textit{Cartesian}. For a binary product $A \times B$, we denote the projection maps $\pi_1^{A,B}$ and $\pi_2^{A,B}$, omitting superscripts when they can be inferred from the context. The pairing of $f : X \rightarrow A$ and $g : X \rightarrow B$ is denoted $\langle f, g \rangle$, and given $\alpha : X \rightarrow A$ and $\beta : Y \rightarrow B$, we write $\alpha \times \beta$ for $\langle \alpha \circ \pi_1, \beta \circ \pi_2\rangle : X \times Y \rightarrow A \times B$. This notation extends to the $n$-fold product $A_1 \times ... \times A_n$, which has projections $\pi_1, ..., \pi_n$ and pairing $\langle f_i \rangle_{i=1}^n$. In particular, for the terminal object $\one$, the unique arrow $X \rightarrow \one$ is denoted $\langle\rangle_X$.

For pullbacks, the notation is similar. The pullback of $f : A \rightarrow C$ and $g : B \rightarrow C$ is denoted $A \times_C B$, and the projections are denoted $\pi_1^{A,B}$ and $\pi_2^{A,B}$. We will also use the notation $\langle p, q \rangle : X \rightarrow A \times_C B$ and $\alpha \times_C \beta : X \times_C Y \rightarrow A \times_C B$ whenever these constructions are well-defined.

\subsubsection*{Functors between slice categories}

Given an arrow $f : A \rightarrow B$ in a category $\mathcal{C}$, we write $\Sigma_f$ for the composition functor $\mathcal{C}/A \rightarrow \mathcal{C}/B$. Its right adjoint is the pullback functor $\Delta_f : \mathcal{C}/B \rightarrow \mathcal{C}/A$ (sometimes denoted $f^*$), assuming it exists. In the special case $B = \one$, these functors become $\Sigma_A : \mathcal{C}/A \rightarrow \mathcal{C}$ (the forgetful functor) and $\Delta_A : \mathcal{C} \rightarrow \mathcal{C}/A$.

If $\Delta_f$ has a further right adjoint, we call it the \textit{dependent product functor} along $f$ and denote it $\Pi_f : \mathcal{C}/A \rightarrow \mathcal{C}/B$. We record the following fact from \cite{Niefield} (corollary 1.2) for later use.

\begin{ffact}
\label{fact:Niefield}
Let $\mathcal{C}$ be a category with finite limits, and let $f : A \rightarrow B$ be an arrow. Then the following are equivalent:
\begin{itemize}
\item $\Delta_f$ has a right adjoint (i.e. $\Pi_f$ exists);
\item $- \times f$ has a right adjoint (i.e. $f$ is exponentiable, per p.8 of \cite{ExpArrow});
\item $\Sigma_A \circ \Delta_f$ has a right adjoint.
\end{itemize}
\end{ffact}

\subsubsection*{Polynomials}

Let $\mathcal{C}$ be a category with finite limits. A \textit{polynomial} over a category $\mathcal{C}$ is a diagram $P$ of the following form.
\begin{equation*}
\begin{tikzcd}
I & A \arrow[l, "s"'] \arrow[r, "f"]
& B \arrow[r, "t"] & J
\end{tikzcd}
\end{equation*}
If the functor $\Pi_f$ exists (i.e. if $f$ is exponentiable, see fact \ref{fact:Niefield}), then the functor $F_P : \mathcal{C}/I \rightarrow \mathcal{C}/J$ given by
\begin{equation*}
\begin{tikzcd}
\mathcal{C}/I \arrow[r, "\Delta_s"]
& \mathcal{C}/A \arrow[r, "\Pi_f"]
& \mathcal{C}/B \arrow[r, "\Sigma_s"]
& \mathcal{C}/J
\end{tikzcd}
\end{equation*}
is called the \textit{polynomial functor associated to $P$}.

A functor $F : \mathcal{C}/I \rightarrow \mathcal{C}/J$ is called a \textit{polynomial functor} if there exists a polynomial $P$ such that $F_P$ exists and $F \cong F_P$. In this case, we say $P$ \textit{represents} $F$. If $F : \mathcal{C} \rightarrow \mathcal{C}$, then we say it is a polynomial functor if this is the case when using the usual identification between $\mathcal{C}$ and $\mathcal{C}/\one$.

These definitions are from \cite{Polynomials}; however, in that paper, they assume that $\mathcal{C}$ is locally Cartesian closed. We avoid this assumption, at the cost of not all polynomial functors existing.

\subsubsection*{Extensivity}

We refer to \cite{Extensivity} for the general definition of an extensive category. Since we will always assume our categories have finite limits, we use the following characterization, drawn from proposition 2.2 of \cite{Extensivity}.

\begin{ffact}
\label{fact:extensivity}
Let $\mathcal{C}$ be a category with finite limits and finite coproducts. Then $\mathcal{C}$ is extensive if and only if for every coproduct $X_1 \longrightarrow X_1 + X_2 \longleftarrow X_2$ and any commutative diagram as below,
\begin{equation*}
\begin{tikzcd}
A_1 \arrow[r] \arrow[d]
& A \arrow[d]
& A_2 \arrow[l] \arrow[d]
\\
X_1 \arrow[r]
& X_1 + X_2 & X_2 \arrow[l]
\end{tikzcd}
\end{equation*}
the squares are pullbacks if and only if the top row is a coproduct.
\end{ffact}

In particular, in an extensive category, every binary coproduct is \textit{universal}, which means that its pullback along any arrow is still a coproduct. Note that whenever we say a category is extensive, we implicitly assume it has finite coproducts.

\section{Reasoning with variables}
\label{sect:IL}

In any Cartesian category, there is a way to express equalities between arrows in a ``variable form", so that an equality $f=g$ is instead written as $f(x) =_x g(x)$. This improves readability: to define some $h : X \times Y \times Z \rightarrow A$, it is much clearer to write
\begin{align*}
h(x,y,z) =_{x,y,z} f(x,g(y,z))
&& \text{than} &&
h = f \circ \langle \pi_1, \; g \circ \langle \pi_2, \pi_3 \rangle \rangle.
\end{align*}
We will use this notation throughout this paper, starting with section \ref{sect:NNOlist}.

The way to make this precise is to develop an \textit{internal language} for our ambient category $\mathcal{C}$, a process which is explained in great detail in sections D1.1 - D1.3 of \cite{Elephant}, particularly on p.837. We will briefly review this process, partly because the language we'll use in this paper (which we take from \cite{Roman}) is even weaker than what is done in \cite{Elephant}, and partly because we also introduce some shorthands for dealing with equalisers and pullbacks.

We should note that, in addition to what is presented in this section, we introduce internal language notation for arithmetic and list-arithmetic operations in sections \ref{sect:NNOlist:ITE}, \ref{sect:NNOlist:order}, and \ref{sect:NNOlist:list}, and we use infix notation when appropriate.

\subsection{Setting up the language}
\label{sect:IL:basics}

First, the internal language has a collection of types, which are the objects of $\mathcal{C}$. It has function symbols of the form $f : A_1 \times ... \times A_n \rightarrow B$, which come from the arrows of $\mathcal{C}$. Arrows with domain $\one$ are treated as constant symbols. We inductively build \textit{typed terms} using variables and these function symbols in the usual way.

Next, a \textit{context} is a list of variables and their types, such as $C=(x:X, y:Y, z:Z)$, and we write $[C]$ for $X \times Y \times Z$. A context is \textit{valid} for a term $t$ if all the variables in $t$ occur in $C$, and the \textit{interpretation} of a term $t:B$ in a valid context $C$ is an arrow $[t]_C : [C] \rightarrow B$. This interpretation is defined inductively: for a variable $x_i$ in a context $C=(x_1:X_1, ..., x_n:X_n)$, $[x_i]_C$ is the projection map $\pi_i : [C] \rightarrow X_i$, and for a term $f(t_1, ..., t_n)$, its interpretation is $f \circ \langle [t_1], ..., [t_n] \rangle$.

Finally, if $t_1, t_2$ are terms in a context $C$, we write $t_1 =_C t_2$ to mean that $[t_1]_C = [t_2]_C$ (as arrows in $\mathcal{C}$). Since this is enough to make sense of the example above with $h(x,y,z)$, this completes the formal setup of the language (note that we don't use relations). Note that the types of the variables in $C$ are omitted if clear from context.

Below, we state some basic facts about this language (the only non-trivial rule is the second one, which requires lemma D1.2.4 in \cite{Elephant}).
\begin{itemize}
\item The equality $=_C$ is an equivalence relation.
\item The equality $=_C$ is preserved by substitution: if $a =_C b$, then $a[\vec{t}/\vec{y}] =_C b[\vec{t}/\vec{y}]$, and if $t_i =_C s_i$ for each $i$, then $a[\vec{t}/\vec{y}] = a[\vec{s}/\vec{y}]$.
\item Given arrows $f,g : X_1 \times ... \times X_n \rightarrow Y$ and the context $\vec{x} = (x_1 : X_1, ..., x_n : X_n)$, we have $f=g$ if and only if $f(x_1, ..., x_n) =_{\vec{x}} g(x_1, ..., x_n)$.
\item For any appropriate arrows $f, (g_i)_{i=1}^m, (h_j)_{j=1}^n$, terms $(t_i)_{i=1}^m$, and context $C$, we have
\begin{equation*}
f\big(g_1(t_1, ..., t_n), ..., g_m(t_1, ..., t_n)\big) =_C (f \circ \langle g_1, ..., g_m \rangle)(t_1, ..., t_n)
\end{equation*}
and
\begin{equation*}
f(h_1(t_1), ..., h_n(t_n)) =_C (f \circ (h_1 \times ... \times h_n))(t_1, ..., t_n).
\end{equation*}
In particular, $f(g(t)) =_C (f \circ g)(t)$.
\end{itemize}

\subsection{Additional notation}
\label{sect:IL:notation}

\subsubsection*{Terms of product types}

Given two terms $a:A$, $b:B$ in a context $C$, we can let $(a,b)$ be shorthand for the term $\Id_{A \times B}(a,b)$ (and similarly for products of any length). Since $(a,b) : A \times B$, this notation lets us mimic products of types without having to build them into the language. Remark that $[(a,b)]_C = \langle [a]_C, [b]_C \rangle$, as we'd expect. In particular, if we have an arrow $g : A \times B \rightarrow X$, then $g((a,b)) =_C g(a,b)$, and so we will not bother to distinguish between these terms.

\subsubsection*{Equalisers}

Suppose we have a context $D = (y_1 : Y_1, ..., y_n : Y_n)$ and terms $s,t : Z$ in the context $D$. Then we can form the following equaliser diagram.
\begin{equation*}
\begin{tikzcd}
{[D']} \arrow[r, phantom, "="]
&[-5ex] E \arrow[r, hookrightarrow, "i"]
& {[D]} \arrow[r, shift left, "{[s]_D}"] \arrow[r, shift right, "{[t]_D}"']
& Z
\end{tikzcd}
\end{equation*}
We will write $\{y_1 : Y_1, ..., y_n : Y_n \mid s=t\}$ as notation for the object $E$, and we will write $(y_1 : Y_1, ..., y_n : Y_n \mid s=t)$ or $(D \mid s=t)$ as notation for the context $D' = (e:E)$. This makes it easier to refer to equalisers without using new variable names; for example, the context for a pullback of $f : U \rightarrow W$, $g : V \rightarrow W$ is $(u:U, \, v:V \mid f(u) = g(v))$.

We want to be talk about ``terms of type $E$" and ``terms of type $[D]$ which equalise $s$ and $t$" interchangeably. This will require us to abuse notation in the following ways. First, if we have a term $e : E$, then we can also think of it as the term $i(e) : [D]$. So, we may write $e : [D]$ to refer to this term.

Second, if we have terms $u_1 : Y_1, ..., u_n : Y_n$ in a context $C = (\vec{x} : \vec{X})$, then the term $(u_1, ..., u_n)$ is formally of type $[D]$. However, if $s[\vec{u}/\vec{y}] =_C t[\vec{u}/\vec{y}]$, we can say $(u_1, ..., u_n)$ is of type $E$ by considering it as shorthand for $\psi(\vec{x}) : E$, where $\psi$ comes from the universal property of the equaliser (see the diagram below).
\begin{equation*}
\begin{tikzcd}
{[D']} \arrow[r, phantom, "="]
&[-5ex] E \arrow[r, hookrightarrow, "i"]
& {[D]} \arrow[r, "{[s]_D}", shift left] \arrow[r, "{[t]_D}"', shift right]
& Z
\\
& {[C]} \arrow[ur, near start, "{\langle [u_1]_C, ..., [u_n]_C \rangle}"']
\arrow[u, "\psi", dashed]
\end{tikzcd}
\end{equation*}
Remark that $[(u_1, ..., u_n)]_C = \psi$.

It is important to note that these two ways of abusing notation are compatible with each other. On one hand, if $(u_1, ..., u_n) : [D]$ equalises $s$ and $t$, then $i(u_1, ..., u_n) =_C (u_1, ..., u_n)$ as terms of type $[D]$; on the other hand, if $e : E$, then $i(e) =_C e$ as terms of type $E$ (these are easy to check).

Finally, we remark that the latter notation behaves well with respect to substitution. The proof is straightforward, though it does require lemma D1.2.4 from \cite{Elephant}.

\begin{rem*}
With the above setup, let $w_1 : X_1, ..., w_m : X_m$ be terms in another context $B$. Then
\begin{align*}
(u_1, ..., u_n)[\vec{w}/\vec{x}] =_B (u_1[\vec{w}/\vec{x}], ..., u_n[\vec{w}/\vec{x}]) \quad (\text{as terms of type } E)
\end{align*}
\end{rem*}

\subsection{Substitution for equalisers}

Let $D = (y_1 : Y_1, ..., y_n : Y_n)$, let $s,t : Z$ be terms in the context $D$, and let $D' = (y_1 : Y_1, ..., y_n : Y_n \mid s=t) = (e:E)$, as in the previous section. Remark that terms in the context $D$ can be also considered to be terms in the context $D'$, simply by ``pulling back" along the inclusion map $i : [D'] \rightarrow [D]$. Notably, if we view $s,t$ as terms in $D'$, then we have $s =_{D'} t$.

Now, consider a term like $g(y_1, ..., y_n)$ in the context $D$. If we consider this term in the context $D'$, can we still substitute terms for the variables $y_i$? The $y_i$ are no longer variables in the context $D'$, so this requires special consideration.

Formally, let $r : R$ be a term in the context $D'$, and let $u_1 : Y_1, ..., u_n : Y_n$ be terms in a context $C$. If we assume that $s[\vec{u}/\vec{y}] =_C t[\vec{u}/\vec{y}]$, then we can say that
\begin{align*}
r[\vec{u}/\vec{y}] \text{ is shorthand for } r[e/\psi(\vec{x})],
\end{align*}
where $\psi$ once again comes from the universal property (see the previous subsection).
This definition ensures that our substitution notation enjoys many of the same properties of regular substitution, and we can check that it behaves the way we want with the following remark (the proof is straightforward).

\begin{rem*}
With the above setup, $y_k[\vec{u}/\vec{y}] =_C u_k$.
\end{rem*}

\section{Basic constructions with natural numbers objects and list objects}
\label{sect:NNOlist}

In this section, we review the standard concepts of \textit{(\param{}) natural numbers objects} and \textit{(\param{}) list objects}. We also describe some basic constructions that can be done with these notions, such as the abstract analogue of the set $E = \{(m,n) \in \N \times \N : m<n\}$, which is needed to state theorem \ref{thm:conclusion}, and establish some basic properties.

\subsection{Natural numbers objects and arithmetic}
\label{sect:NNOlist:NNO}

We start this section by recalling the standard definition of natural numbers objects, which can be found in \cite{Roman}, for one.

\begin{dfn}
Let $\mathcal{C}$ be a Cartesian category. A \textit{(\param{}) natural numbers object} (NNO) is a triple $(N,0,s)$ consisting of an object $N$ and arrows $0 : \one \rightarrow N$, $s : N \rightarrow N$ which satisfy the following property. For any objects $A,B$ and arrows $g : A \rightarrow B$, $h : A \times N \times B \rightarrow B$, there exists a unique arrow $f : A \times N \rightarrow B$ such that $f(a,0) =_a g(a)$ and $f(a,sn) =_{a,n} h(a,n,f(a,n))$.
\end{dfn}

This definition uses our internal language (see section \ref{sect:IL}); the definition is given diagrammatically in \cite{Roman}. Following the convention of \cite{Roman}, the term ``natural numbers object" will always refer to \textit{\param{}} natural numbers objects; we will not discuss un\param{} NNOs (where we omit the parameter $A$ in the definition).

The NNO property lets us define the operations of addition $+ : N \times N \rightarrow N$, multiplication $\bullet : N \times N \rightarrow N$, predecessor $P : N \rightarrow N$, and (truncated) subtraction $\dotminus : N \times N \rightarrow N$ in the usual inductive way. The operations $+$ and $\bullet$, along with $0$ and $1 = s0$, make $N$ into an (internal) commutative semiring. More properties of these operations can be found in \cite{Roman} (themselves inspired by \cite{Goodstein}), and we prove some additional facts in appendix \ref{app:NNO}.

Finally, it's helpful to have the following formulation of the universal property that uses the internal language even more (we omit the proof).
\begin{rem*}
Let $N$ be a natural numbers object. Let $A,B$ be objects, let $g$ be a term in the context $(a:A)$, and let $h$ be a term in the context $(a:A,n:N,b:B)$. Then there exists a term $f$ in the context $(a:A,n:N)$ such that $f[0/n] =_a g$ and $f[s(n)/n] =_{a,n} h[f/b]$. Moreover, if $f_1, f_2$ are two such terms, then $f_1 =_{a,n} f_2$.
\end{rem*}
This avoids the need to think of the terms $f, f_1, f_2, g, h$ as coming from arrows, and is more true to how calculations are performed in \cite{Roman}.

\subsection{Definition by cases in NNOs}
\label{sect:NNOlist:ITE}

In this section, we'll establish some notation for defining functions by cases. As motivation, we start with a simple fact about natural numbers objects in Cartesian categories (see, e.g., \cite[lemma A2.5.5]{Elephant}).

\begin{ffact}
\label{fact:coprodNNO}
In a Cartesian category with a \param{} NNO, the diagram $\one \arrowset{0} N \arrowsetl{s} N$ is a coproduct.
\end{ffact}

This fact suggests that we can use an NNO to define functions by cases. We start by taking an object $B$ and defining the function $ITE_B : B \times B \times N \rightarrow B$ (short for ``if-then-else") by induction as follows:
\begin{align*}
ITE_B(x,y,0) =_{x,y} x,
&& ITE_B(x,y,sn) =_{x,y,n} y.
\end{align*}
We will use the following special notation for this function.
\begin{align*}
ITE_B(x,y,n) =_{x,y,n} \left\{\begin{matrix}
x & \text{if } n=0 \\
y & \text{else}
\end{matrix}\right.
\end{align*}
For legibility, we will introduce some notation for when we want to consider more than one case at a time. We'll write
\begin{align*}
\left\{\begin{matrix}
x & \text{if } m=0 \\
y & \text{else if } n=0 \\
z & \text{else}
\end{matrix}\right\}
=_{x,y,z,m,n}
\left\{\begin{matrix}
x & \text{if } m=0 \\
\left\{\begin{matrix}
y & \text{if } n=0 \\
z & \text{else}
\end{matrix}\right\} & \text{else}
\end{matrix}\right\} .
\end{align*}
Finally, in some cases, we will want our condition to look like $m \leq n$ or $m<n$ instead of $n = 0$; we'll describe how this works in section \ref{sect:NNOlist:order}.

\subsection{Notation for order in NNOs}
\label{sect:NNOlist:order}

An important facet of arithmetic is order, and in this section we introduce notation to talk about it in categories with natural numbers objects.

Given two natural numbers $m,n$, we know that $m \leq n$ if and only if $m \dotminus n = 0$. Therefore, given two terms $m,n : N$ in a context $C$, we write $m \leq_C n$ to mean that $m \dotminus n =_C 0$, and we write $m <_C n$ to mean $s(m) \dotminus n =_C 0$. This lets us write things like $x <_x x+1$.

We will require various facts about these order relations. Most are left to appendix \ref{app:NNO}, but we include the following proposition here because it will be used several times.

\begin{fprop}
\label{prop:>0}
Let $\mathcal{C}$ be a Cartesian category with NNO, and let $n : N$ be a term in a context $C$. If $n >_C 0$, then $n =_C s(P(n))$.
\end{fprop}

\begin{sproof}
If $n >_C 0$, then $n \geq_C s(0) = 1$, so $n =_C m+1 =_C s(m)$ for some term $m : N$ (by proposition \ref{prop:appN:leqEquiv} in appendix \ref{app:NNO}). Thus $s(P(n)) =_C s(P(s(m))) =_C s(m) =_C n$ (since $P(s(m)) =_C m$ by definition of $P$).
\end{sproof}

There are two more cases where we may write $m \leq n$ and $m < n$ as shorthand for $m \dotminus n = 0$ and $sm \dotminus n = 0$, respectively. The first is when we denote equalisers in the internal language: we will write, for instance, $\{m:N, n:N \mid m < n\}$ instead of $\{m:N, n:N \mid (sm \dotminus n) = 0\}$. The second is when using the notation of section \ref{sect:NNOlist:ITE}; we add the following shorthand (and similarly for $m<n$).
\begin{align*}
\left\{\begin{matrix}
x & \text{if } m \leq n \\
y & \text{else}
\end{matrix}\right\}
=_{x,y,m,n}
\left\{\begin{matrix}
x & \text{if } (m \dotminus n)=0 \\
y & \text{else}
\end{matrix}\right\}
\end{align*}

\subsection{The collection of finite cardinals}
\label{sect:NNOlist:E}

A key part of theorem \ref{thm:conclusion} is the arrow $\pi_2^E : E \rightarrow N$. In the category of sets, we said that $E = \{(m,n) \in \N \times \N : m < n\}$, and this carries over nicely to a category with finite limits and a natural numbers object $N$, where we set
\begin{align*}
E = \{m : N, \, n : N \mid m<n \}.
\end{align*}
(Here, we use notation from sections \ref{sect:IL:notation} and \ref{sect:NNOlist:order}; explicitly, $E$ is the equalizer of the arrows $N \times N \rightarrow N$ given by $(m,n) \mapsto s(m) \dotminus n$ and $(m,n) \mapsto 0$.) Now, if we write $D$ for the context $(m : N, \, n : N \mid m<n)$ (remark that $[D] = E$), then we can define the map $\pi_2^E : E \rightarrow N$ by $\pi_2^E(m,n) =_D n$. Going forward, we'll use the notation $E, \pi_2^E$ to refer to this object without mentioning it explicitly.

\begin{rem*}
It is not hard to convince yourself that the above equalizer exists even if we only assume the existence of finite products (instead of finite limits): $E$ is a countably infinite set, and is therefore isomorphic to $N$. However, checking this formally is technical and not relevant to the rest of this paper, so we omit the proof of this.
\end{rem*}

\pagebreak

\subsection{List objects}
\label{sect:NNOlist:list}

In this section, we introduce list objects and give some basic properties. List objects in the context of category theory originated with Cockett in, for instance, \cite{CockettList}; for our definition, we use the version given in \cite{MaiettiList}. As with natural numbers objects, ``list object" will always mean the \textit{\param{}} version.

\begin{dfn}
Let $\mathcal{C}$ be a Cartesian category, and let $X \in \mathcal{C}$ be an object. A \textit{(\param{}) list object} on $X$ is a triple $(L(X), r_0^X, r_1^X)$ consisting of an object $L(X)$ and arrows $r_0^X : \one \rightarrow L(X)$, $r_1^X : X \times L(X) \rightarrow L(X)$ which satisfy the following property.
\begin{quote}
For any objects $A,B$ and arrows $g : A \rightarrow B$, $h : X \times B \rightarrow B$, there exists a unique arrow $f : A \times L(X) \rightarrow B$ which makes the following diagrams commute.
\begin{equation*}
\begin{tikzcd}
&[4ex] A \times L(X) \arrow[dd, "f"] \\[-2ex]
A \arrow[ur, "{\langle \Id_A, r_0^X \rangle}"]
\arrow[dr, "g"'] & \\[-2ex]
& B
\end{tikzcd}
\qquad \qquad
\begin{tikzcd}
A \times X \times L(X) \arrow[r, "\Id_A \times r_1^X"]
\arrow[d, "{\langle \pi_X, f \rangle}"']
&[5ex] A \times L(X) \arrow[d, "f"]
\\[5ex]
X \times B \arrow[r, "h"']
& B
\end{tikzcd}
\end{equation*}
\end{quote}
We ignore the slight abuses of notation in these diagrams because we prefer to write them equationally. By writing $\varnothing$ for $r_0^X$ and $x :: \ell$ for $r_1^X(x,\ell)$, the diagrams become the equations $f(a,\varnothing) =_a g(a)$ and $f(a, x :: \ell) =_{a,x,\ell} h(x,f(a,\ell))$. A version of the universal property with terms can also be formulated, like for natural numbers objects (see section \ref{sect:NNOlist:NNO}).
\end{dfn}

The list object property can actually be generalized slightly to the following proposition.

\begin{fprop}
Let $\mathcal{C}$ be a Cartesian category which has all list objects, and let $X \in \mathcal{C}$. For any objects $A,B$ and arrows $g : A \rightarrow B$, $h : A \times X \times L(X) \times B \rightarrow B$, there exists a unique arrow $f : A \times L(X) \rightarrow B$ such that $f(a, \varnothing) =_a g(a)$ and $f(a, x :: \ell) =_{a,x,\ell} h(a,x,\ell,f(a,\ell))$.
\end{fprop}

\begin{sproof}
In the definition of list objects, replace $B$ with $A \times L(X) \times B$. With some small straightforward adjustments, we get the desired property.
\end{sproof}

From this proposition, it becomes clear that $L(\one)$ is a natural numbers object. This means that it is redundant to assume that a category has list objects \textit{and} a natural numbers object.

Similarly to natural numbers objects, list objects can be expressed as a coproduct. While this fact can be found in \cite{CockettList}, it is stated in a much more abstract way, so we provide a simple proof.

\begin{ffact}
\label{fact:coprodList}
In a Cartesian category with list objects, the diagram $\one \arrowset{r_0^X} L(X) \arrowsetl{r_1^X} X \times L(X)$ is a coproduct.
\end{ffact}

\begin{sproof}
With the generalized version of the list object property, if we set $A = \one$, $g = u_0$, and $h = u_1 \circ \langle \pi_2, \pi_3 \rangle$, the diagrams become the following.
\begin{equation*}
\begin{tikzcd}
\one \arrow[r, "r_0^X"] \arrow[dr, "u_0"']
& L(X) \arrow[d, dashed, "f'"]
& X \times L(X) \arrow[l, "r_1^X"'] \arrow[dl, "u_1"]
\\
& B
\end{tikzcd}
\end{equation*}
This is the usual coproduct diagram.
\end{sproof}

Finally, if a Cartesian category $\mathcal{C}$ has all list objects, then we can form a functor $L : \mathcal{C} \rightarrow \mathcal{C}$ which maps $X$ to the list object $L(X)$. Its action on arrows is defined as follows: given $f : X \rightarrow Y$, we define $L(f) : L(X) \rightarrow L(Y)$ inductively by
\begin{align*}
L(f)(\varnothing) = \varnothing
&& \text{and} &&
L(f)(x :: \ell) =_{x,\ell} f(x) :: L(f)(\ell).
\end{align*}

\subsection{Important list arrows}
\label{sect:NNOlist:listars}

In this section, we define some important arrows related to list objects and prove some properties about them. Throughout, we assume we are working in a Cartesian category $\mathcal{C}$ with \param{} list objects (and therefore a \param{} NNO).

We start with the definitions. These operators are all indexed by an object $X$, but we remove this index if it is clear from context.

\begin{itemize}
\item The \textit{length} arrow, $\len_X : L(X) \rightarrow N$, is defined inductively by 
\begin{align*}
\len_X(\varnothing) = 0, \qquad
\len_X(x :: \ell) =_{x,\ell} s(\len_X(\ell)).
\end{align*}

\item The \textit{truncate} arrow, $\tr_X : L(X) \rightarrow L(X)$, is the arrow which removes the leading element from a list. It is defined inductively by
\begin{align*}
\tr_X(\varnothing) = \varnothing,
\qquad
\tr_X(x :: \ell) =_{x,\ell} \ell.
\end{align*}

\item The \textit{tail} arrow, $\tail_X : N \times L(X) \rightarrow L(X)$, is an arrow which iterates truncation: $\tail_X(n,\ell)$ returns the list $\ell$ with the $n$ leading elements removed. It is defined inductively by
\begin{align*}
\tail_X(0,\ell) =_\ell \ell,
\qquad
\tail_X(s(n), \ell) =_{n,\ell} \tr_X(\tail_X(n,\ell)).
\end{align*} 

\item The \textit{zeroth-or-default} arrow, $\zerothDef_X : X \times L(X) \rightarrow X$, is an arrow which takes a pair $(x,\ell)$ as input. If $\ell$ is a non-empty list, it returns its leading element; however, if $\ell$ is empty, it returns $x$, the ``default" element. This arrow is defined inductively by
\begin{align*}
\zerothDef_X(x,\varnothing) =_x x,
\qquad
\zerothDef_X(x, y :: \ell) =_{x,y,\ell} y.
\end{align*}
\end{itemize}
Note that we say ``the zeroth element" for the leading element of a list instead of ``the first element". This is to avoid confusion later when we say ``the $n^{th}$ element", since ``the $n^{th}$ element" with $n=0$ will refer to the leading element.

\begin{fprop}
\label{prop:Lfunctor}
Let $f : X \rightarrow Y$ be an arrow. Then
\begin{itemize}
\item $\len_Y\big(L(f)(\ell)\big) =_\ell \len_X(\ell)$;
\item $\tr_Y\big(L(f)(\ell)\big) =_\ell L(f)\big(\tr_X(\ell)\big)$;
\item $\tail_Y\big(m, \, L(f)(\ell)\big) =_{m,\ell} L(f)\big(\tail_X(m, \, \ell)\big)$;
\item $\zerothDef_Y\big(f(x), \, L(f)(\ell)\big) =_{x,\ell} f\big(\zerothDef_X(x, \, \ell)\big)$.
\end{itemize}
\end{fprop}

\begin{sproof}
Each of these equalities is proved by list induction. For the first, we show $\len_Y \circ L(f)$ satisfies the equations defining $\len_X$. For the third, we show that both terms satisfy $h(0,\ell) =_\ell L(f)(\ell)$ and $h(s(n),\ell) =_{n,\ell} \tr_Y(h(n,\ell))$, using the previous result about $\tr$. The other two are simple.
\end{sproof}

The next proposition shows how to decompose list objects based on their length. We draw particular attention to the second case: we would like to decompose a list $t : L(X)$ of positive length into $t = x :: t'$, but this is not possible unless we assume \ref{hyp:list} (see proposition \ref{prop:H1term}). The key limitation is that we cannot extract a term of type $X$ from the context $(\ell : L(X) \mid \len(\ell) > 0)$ without additional assumptions (see also proposition \ref{prop:def}), so in order to perform the decomposition, we need to be given an arbitrary term $\deff : X$.

\begin{fprop}
\label{prop:listLenDecomp}
Let $t : L(X)$ be a term in a context $C$.
\begin{enumerate}
\item If $\len(t) =_C 0$, then $t =_C \varnothing$.
\item If $\len(t) >_C 0$, then for any term $\deff : X$ in the context $C$, we have
\begin{align*}
t =_C \zerothDef(\deff, t) :: \tr(t).
\end{align*}
\end{enumerate}
\end{fprop}

\begin{sproof}
For the first point, consider the map $Z : N \times L(X) \rightarrow L(X)$ defined inductively by $Z(0,\ell) =_{\ell} \varnothing$, $Z(sn,\ell) =_{n,\ell} \ell$. By induction on $\ell$, it is easy to show that $Z(\len(\ell), \ell) =_\ell \ell$. We conclude by using the hypothesis $\len(t) =_C 0$ to compute
\begin{align*}
t =_C Z(\len(t), t) =_C Z(0,t) =_C \varnothing,
\end{align*}
as desired.

For the second point, we similarly consider an arrow $Z : X \times N \times L(X) \rightarrow L(X)$ defined by $Z(x,0,\ell) =_{x,\ell} \varnothing$ and $Z(x,sn,\ell) =_{x,n,\ell} \zerothDef(x,\ell) :: \tr(\ell)$. As before, we can easily show by induction that $Z(x, \len(\ell), \ell) =_{x,\ell} \ell$. Now, if $\len(t) >_C 0$, then $\len(t) =_C s(n)$ for some $n : N$ (by proposition \ref{prop:>0}). If moreover we have some $\deff : X$, then
\begin{align*}
t =_C Z(\deff, \len(t), t)
=_C Z(\deff, s(n), t)
=_C \zerothDef(\deff,t) :: \tr(t),
\end{align*}
as desired.
\end{sproof}

Finally, we check that the $\tail$ operation has the length we expect. Combined with the previous result, this gives us a way to tell when $\tail(n,\ell)$ is the empty list.

\begin{fprop}
\label{prop:lenIterTr}
We have $\len(\tail(n,\ell)) =_{n,\ell} \len(\ell) \dotminus n$. In particular, for any terms $n : N$ and $\ell : L(X)$ in a context $C$, if $n \geq_C \len(\ell)$, then $\tail(n,\ell) =_C \varnothing$.
\end{fprop}

\begin{sproof}
First, we note that $\len(\tr(\ell)) =_\ell P(\len(\ell))$; this is easily checked by induction on $\ell$.
Using this fact, it is easy to check that the terms $\len(\tail(n,\ell))$ and $\len(\ell) \dotminus n$ both satisfy $h(0,\ell) =_\ell \len(\ell)$ and $h(sn,\ell) =_{n,\ell} P(h(n,\ell))$, so we get the desired equality by induction.

For the second part, let $n : N$ and $\ell : L(X)$ be terms in a context $C$. If $n \geq_C \len(\ell)$, then $\len(\tail(n,\ell)) =_C \len(\ell) \dotminus n =_C 0$, so $\tail(n,\ell) =_C \varnothing$ by proposition \ref{prop:listLenDecomp}.
\end{sproof}

\section{Proof setup and outline}
\label{sect:outline}

Recall that the goal of this paper is to prove theorem \ref{thm:conclusion}, which states that if $\mathcal{C}$ is an extensive category with finite limits and list objects, then the list object functor $L$ is polynomial.

We start by proving the following proposition, which makes this goal more attainable. We make use of the functor $L_N : \mathcal{C} \rightarrow \mathcal{C}/N$ which maps an object $X$ to the arrow $\len_X : L(X) \rightarrow N$.

\begin{fprop}
\label{prop:strategy}
Let $\mathcal{C}$ be a category with finite limits and \param{} list objects. To show that the polynomial $\one \leftarrow E \arrowset{\pi_2^E} N \rightarrow \one$ represents $L$, it suffices to show that $L_N$ is a right adjoint to $\Sigma_E \circ \Delta_{\pi_2^E} : \mathcal{C}/N \rightarrow \mathcal{C}$.
\end{fprop}

\begin{sproof}
Suppose $L_N$ is right adjoint to $\Sigma_E \circ \Delta_{\pi_2^E}$. By fact \ref{fact:Niefield} (due to \cite{Niefield}), this right adjoint existing implies that $\Pi_{\pi_2^E}$, the right adjoint of $\Delta_{\pi_2^E}$, also exists. Thus the polynomial functor associated to $\one \leftarrow E \arrowset{\pi_2^E} N \rightarrow \one$ exists, and it remains to show that
\begin{align*}
\Sigma_N \circ \Pi_{\pi_2^E} \circ \Delta_E \cong L.
\end{align*}
We note that $L = \Sigma_N \circ L_N$, so we just need to show that $\Pi_{\pi_2^E} \circ \Delta_E \cong L_N$. For this, we note that $\Pi_{\pi_2^E}$ and $\Delta_E$ are the right adjoints of $\Delta_{\pi_2^E}$ and $\Sigma_E$, respectively, and so $\Pi_{\pi_2^E} \circ \Delta_E$ is right adjoint to $\Sigma_E \circ \Delta_{\pi_2^E}$ (see \cite[prop. 3.2.1]{Borceux1}). However, $L_N$ is already a right adjoint of $\Sigma_E \circ \Delta_{\pi_2^E}$ by assumption, so by uniqueness of adjoints, we get $\Pi_{\pi_2^E} \circ \Delta_E \cong L_N$, as desired.
\end{sproof}

Thus, our proof strategy will involve showing that $L_N$ is a right adjoint to $\Sigma_E \circ \Delta_{\pi_2^E}$. Note that the above proposition does not require $\mathcal{C}$ to be extensive; in fact, we can weaken the assumption of extensivity and instead use the following two hypotheses. (In section \ref{sect:technical}, we rephrase these hypotheses using the internal language to make them easier to use.)

\begin{enumerate}
\item[\mylabel{hyp:list}{(H1)}]
For any object $X$, in the following commutative diagram, the squares are pullbacks.
\begin{equation*}
\begin{tikzcd}
\one \arrow[r, "r_0^X"] \arrow[d]
& L(X) \arrow[d, "\len_X"]
& X \times L(X) \arrow[l, "r_1^X"'] \arrow[d, "\len_X \circ \pi_2"]
\\
\one \arrow[r, "0"]
& N & N \arrow[l, "s"']
\end{tikzcd}
\end{equation*}

\item[\mylabel{hyp:NNO}{(H2)}]
The coproduct $\begin{tikzcd} \one \arrow[r,"0"] & N & N \arrow[l,"s"'] \end{tikzcd}$ is universal.
\end{enumerate}

\begin{fprop}
\label{prop:extImpliesH1H2}
Let $\mathcal{C}$ be a category with finite limits and \param{} list objects. If $\mathcal{C}$ is extensive, then hypotheses \ref{hyp:list} and \ref{hyp:NNO} hold.
\end{fprop}

\begin{sproof}
This is a direct consequence of the characterization of extensivity mentioned in section \ref{sect:prelim} (i.e. fact \ref{fact:extensivity}). For \ref{hyp:list}, the squares are pullbacks because the top and bottom rows are coproducts (facts \ref{fact:coprodNNO} and \ref{fact:coprodList}) and the diagram commutes (by definition of $\len$).
\end{sproof}

Note that hypothesis \ref{hyp:NNO} is not trivial: the coproduct $\one \arrowset{0} N \arrowsetl{s} N$ is not necessarily universal in categories with finite limits. A counter-example can be found on the MathOverflow website \cite{MOcounterex}. On the other hand, in hypothesis \ref{hyp:list}, the left square can be easily shown to be a pullback without extra assumptions, so assuming it is redundant. We only need to assume that the right square is a pullback; though we have not constructed a counter-example, we suspect this assumption is necessary.

With this setup, the rest of the paper will be devoted to proving the following theorem.

\begin{fthm}
\label{thm:finalPolyn}
Let $\mathcal{C}$ be a category with finite limits and \param{} list objects. If we assume \ref{hyp:list} and \ref{hyp:NNO}, then $L_N$ is a right adjoint to $\Sigma_E \circ \Delta_{\pi_2^E}$.
\end{fthm}

Our strategy for this theorem is to use the universal property characterization and construct a natural transformation $\nth : \Sigma_E \circ \Delta_{\pi_2^E} \circ L_N \Rightarrow \Id_{\mathcal{C}}$ such that $(L_N(X), \nth_X)$ is a universal morphism from $\Sigma_E \circ \Delta_{\pi_2^E}$ to $X$ for each $X$. Unravelling the definitions, this means that we must prove the following universal property for each $X$:
\begin{center}
For any $l_A : A \rightarrow N$ and $g : E \times_N A \rightarrow X$, there exists a unique $h : A \rightarrow L(X)$ such that $l_A = \len_X \circ h$ and $g = \nth_X \circ (\Id \times_N h)$.
\end{center}
So, we can think of the proof as being broken down into three steps.
\begin{enumerate}
\item First, we must define the arrows $\nth_X : E \times_N L(X) \rightarrow X$ and check that they form a natural transformation.

\item Second, we must construct $h$ based on the given $l_A$ and $g$, and check that it satisfies the required equations.

\item Third, we must show that any two arrows satisfying the equations must be equal.
\end{enumerate}
These steps are accomplished mostly independently: constructing $\nth_X$ is done in section \ref{sect:nth}; constructing the map $h$ is done in section \ref{sect:list}; and showing uniqueness is done in section \ref{sect:equality}. Section \ref{sect:technical} provides some important technical tools for these steps. Putting it all together gives us the following proof.

\pagebreak

\begin{sproof}[Proof of theorem \ref{thm:finalPolyn}]
Proposition \ref{prop:epsNT} tells us that $\nth : \Sigma_E \circ \Delta_{\pi_2^E} \circ L_N \Rightarrow \Id$ is a natural transformation. Moreover, theorem \ref{thm:UPexistence} and corollary \ref{corol:listEqExt} together tell us that for each $X$, $(L_N(X), \nth_X)$ is a universal morphism from $\Sigma_E \circ \Delta_{\pi_2^E}$ to $X$. Thus, $L_N$ is right adjoint to $\Sigma_E \circ \Delta_{\pi_2^E}$ (see e.g. \cite[dfn. 3.1.4]{Borceux1}).
\end{sproof}

Of course, this now implies theorem \ref{thm:conclusion}.

\begin{sproof}[Proof of theorem \ref{thm:conclusion}]
By assumption, we're working in an extensive category, so \ref{hyp:list} and \ref{hyp:NNO} hold by remark \ref{prop:extImpliesH1H2}. We can thus apply theorem \ref{thm:finalPolyn} to conclude that $L_N$ is a right adjoint to $\Sigma_E \circ \Delta_{\pi_2^E}$, and proposition \ref{prop:strategy} tells us that this implies the desired conclusion.
\end{sproof}

\section{Technical calculations}
\label{sect:technical}

In this section, we prove some technical results that will be used at various points in the rest of the paper. Throughout, we assume that we are working in a category $\mathcal{C}$ with finite limits and \param{} list objects (and thus a \param{} natural numbers object).

\subsection{Rephrasing the hypotheses}
\label{sect:technical:H1H2}

Hypotheses \ref{hyp:list} and \ref{hyp:NNO} are not easy to apply directly when doing proofs with the internal language. In this section, we obtain consequences of these hypotheses that can be applied to the internal language. We start with \ref{hyp:list}, which is fairly easy to adapt.

\begin{fprop}
\label{prop:H1term}
Assume \ref{hyp:list}. If $\ell : L(X)$ is a term in a context $C$ such that $\len(\ell) >_C 0$, then there exist terms $x : X$, $\ell' : L(X)$ in the context $C$ such that $\ell =_C x :: \ell'$.
\end{fprop}

\begin{sproof}
Since $\len(\ell) >_C 0$, we have $\len(\ell) =_C s(P(\len(\ell)))$ by proposition \ref{prop:>0}, which means the following diagram commutes.
\begin{equation*}
\begin{tikzcd}
& & {[C]}
\arrow[dll, "{[\ell]_C}"', bend right]
\arrow[ddl, "{P \circ \len \circ [\ell]_C}", bend left]
\\
L(X) \arrow[d, "\len_X"']
& X \times L(X) \arrow[d] \arrow[l, "r_1^X"']
\\
N & N \arrow[l, "s"']
\end{tikzcd}
\end{equation*}
Since \ref{hyp:list} asserts that this square is a pullback, there exists an arrow $\gamma : [C] \rightarrow X \times L(X)$ which makes the whole diagram commute. So, if we write $C  = (y_1 : Y_1, ..., y_n : Y_n)$ and we let $x =_C \pi_1(\gamma(y_1, ..., y_n))$ and $\ell' =_C \pi_2(\gamma(y_1, ..., y_n))$, then the equality $[\ell]_C = r_1^X \circ \gamma$ implies $\ell =_C x :: \ell'$, as desired.
\end{sproof}

Next, we move on to \ref{hyp:NNO}. This hypothesis is used when we need to perform a construction or calculation ``by cases", depending on whether a certain parameter is zero or non-zero (e.g. the proof of theorem \ref{thm:UPexistence}); the following proposition makes this easier to do.

\begin{fprop}
\label{prop:coproduct}
Let $t : N$ be a term in a context $C = (x_1 : X_1, ..., x_n : X_n)$, and let $C_0$, $C_{>0}$ be the contexts $C_0 = (C \mid t=0)$ and $C_{>0} = (C \mid t > 0)$.

If we assume hypothesis \ref{hyp:NNO}, then
$\begin{tikzcd}
{[C_0]} \arrow[r, hookrightarrow]
& {[C]} \arrow[r, hookleftarrow]
& {[C_{>0}]} 
\end{tikzcd}$
is a coproduct.
\end{fprop}

\begin{sproof}
Write $i_0, i_{>0}$ for the inclusion maps of $[C_0], [C_{>0}]$ into $[C]$, respectively, and consider the following diagram.
\begin{equation*}
\begin{tikzcd}
{[C_0]} \arrow[r, "i_0"] \arrow[d]
& {[C]} \arrow[d, "{[t]_C}"]
& {[C_{>0}]} \arrow[l, "i_{>0}"'] \arrow[d, "{P \circ [t]_C \circ i_{>0}}"]
\\
\one \arrow[r, "0"]
& N & N \arrow[l, "s"']
\end{tikzcd}
\end{equation*}
We claim that this diagram commutes and that the squares are pullbacks. If this is true, then hypothesis \ref{hyp:NNO} tells us that the top row is a coproduct, as desired.

First, we check that the squares commute. The left square states that, in the context $C_0$, $t =_{C_0} 0$, which is true by definition. The right square states that, in the context $C_{>0}$, $t =_{C_{>0}} s(P(t))$. We know $t > 0$ in this context, so the equality follows by proposition \ref{prop:>0}.

Second, we claim that the left square is a pullback because of how $[C_0]$ is defined as an equalizer, and the right square is a pullback because of how $[C_{>0}]$ is defined as an equalizer. This is mostly a routine check, so we omit the details. We just note that for the latter case, if $f : B \rightarrow [C]$ and $g : B \rightarrow N$ are such that $[t]_C \circ f = s \circ g$, then $[t]_C \circ f > 0$ (by proposition \ref{prop:>0}) and $g = P \circ [t]_C \circ f$.
\end{sproof}

The conclusion of this proposition can be reformulated using the internal language so that it is easier to use.

\begin{rem*}
Let $t, C, C_0, C_{>0}$ be as in proposition \ref{prop:coproduct}, and assume \ref{hyp:NNO}.
Let $a_0 : Y$ be a term in the context $C_0$, and let $a_{>0}$ be a term in the context $C_{>0}$. Then there exists a unique arrow $h : [C] \rightarrow Y$ such that $h(x_1, ..., x_n) =_{C_0} a_0$ and $h(x_1, ..., x_n) =_{C_{>0}} a_{>0}$.
Moreover, for any terms $p, q : Y$ in the context $C$, if $p =_{C_0} q$ and $p =_{C_{>0}} q$, then $p =_C q$.
\end{rem*}

We will generalise the above result slightly. Instead of splitting into the case $t=0$ and $t>0$, we will split into the case $u<w$ and $u \geq w$. Note that the conclusion of this generalization can also be phrased using terms as in the above remark; we won't rewrite it in full.

\begin{fcorol}
\label{corol:coproductLeq}
Let $u,w : N$ be terms in a context $C = (x_1 : X_1, ..., x_n : X_n)$, and let $C_1$, $C_2$ be the contexts $C_1 = (C \mid u<w)$ and $C_2 = (C \mid u \geq w)$.

If we assume hypothesis \ref{hyp:NNO}, then
$\begin{tikzcd}
{[C_1]} \arrow[r, hookrightarrow]
& {[C]} \arrow[r, hookleftarrow] & {[C_2]}
\end{tikzcd}$
is a coproduct.
\end{fcorol}

\begin{sproof}
Remark that $C_2 = (C \mid w \dotminus u = 0)$. So, if we let $C_{>0} = (C \mid w \dotminus u > 0)$, then by proposition \ref{prop:coproduct}, the following diagram is a coproduct.
\begin{equation*}
\begin{tikzcd}
{[C_2]} \arrow[r, hookrightarrow]
& {[C]} \arrow[r, hookleftarrow]
& {[C_{>0}]}
\end{tikzcd}
\end{equation*}
So, to show the result we want, we need to show that there is an isomorphism $[C_1] \cong [C_{>0}]$ which respects the inclusion into $[C]$.

The only choice for an isomorphism between the objects $[C_1] = \{C \mid u < w\}$ and $[C_{>0}] = \{C \mid w \dotminus u > 0\}$ which respects the inclusion into $[C]$ is the ``identity map" $(x_1, ..., x_n) \mapsto (x_1, ..., x_n)$. Therefore, we just need to show that this map is well-defined in both directions. This just involves showing that $w \dotminus u >_{C_1} 0$ (knowing that $u <_{C_1} w$) and that $u <_{C_{>0}} w$ (knowing that $w \dotminus u >_{C_{>0}} 0$). This is true by proposition \ref{prop:appN:calculation} of the appendix.
\end{sproof}

\subsection{Pullbacks along E}
\label{sect:technical:contextE}

As we saw in proposition \ref{prop:strategy}, a central part of our proof strategy is the functor $\Sigma_E \circ \Delta_{\pi_2^E} : \mathcal{C}/N \rightarrow \mathcal{C}$, which maps an arrow $l_A : A \rightarrow N$ to the object $E \times_N A$ defined by the following pullback.
\begin{equation*}
\begin{tikzcd}
E \times_N A \arrow[r] \arrow[d] & A \arrow[d, "l_A"]
\\ E \arrow[r, "\pi_2^E"] & N
\end{tikzcd}
\end{equation*}
Using the internal language, we have $E \times_N A = \{e : E, \, a : A \mid \pi_2^E(e) = l_A(a) \}$. Since $E$ is the object $\{m : N, n : N \mid m < n\}$, we could also write this as $\{(m,n) : E, \, a : A \mid n = l_A(a)\}$. However, this shows us that the ``$n$" parameter is redundant, and that we can just think of this object as $\{m : N, \, a : A \mid m < l_A(a)\}$. Specifically, we'd like to use the context $(m:N, a:A \mid m < l_A(a))$ when referencing the object $E \times_N A$. The following proposition assures us that this is valid.

\begin{fprop}
\label{prop:Eeq}
Let $l_A : A \rightarrow N$. Then the following maps are well-defined and provide an isomorphism $E \times_N A \cong \{m : N, \, a : A \mid m < l_A(a)\}$.
\begin{equation*}
\begin{tikzcd}
E \times_N A
\arrow[r, bend left, "{((m,n),a) \mapsto (m,a)}", start anchor = north east, end anchor = north west]
&[6ex]
\{m : N, \, a : A \mid m < l_A(a)\}
\arrow[l, bend left, "{((m,l_A(a)),a) \mapsfrom (m,a)}", start anchor = south west, end anchor = south east]
\end{tikzcd}
\end{equation*}
\end{fprop}

\begin{sproof}
Straightforward check.
\end{sproof}

Note that we do need to be careful when writing $\{m : N, a : A \mid m < l_A(a)\}$ for $E \times_N A$, because applying the projection maps does not behave the way we expect. For instance, we have the map $\pi_1^{E,A} : E \times_N A \rightarrow E$, and if we consider $(m,a) : E \times_N A$ with $m : N, a : A$, then $\pi_1^{E,A}(m,a)$ is not equal to $m$, as we might think if we're not being careful.

Fortunately, this issue does not come up very often. We'll only need the following fact to assure us that we aren't making mistakes (the proof is a straightforward check, so we omit it).

\begin{ffact}
\label{fact:Idxf}
Suppose we have arrows $l_A : A \rightarrow N$, $l_B : B \rightarrow N$ and $f : A \rightarrow B$. Assume $l_B \circ f = l_A$, so that we can consider the following arrow.
\begin{equation*}
\begin{tikzcd}
E \times_N A \arrow[r, "\Id_E \times_N f"]
&[6ex] E \times_N B
\end{tikzcd}
\end{equation*}
Then, in the context $C = (m : N, a : A \mid m < l_A(a))$, we have
\begin{align*}
(\Id_E \times_N f)(m,a) =_C (m,f(a)).
\end{align*}
\end{ffact}

\subsection{Converting a finite sequence to an infinite sequence}
\label{sect:technical:finiteToInf}

In the proof outline of section \ref{sect:outline}, we saw that we need to use an arrow $g : E \times_N A \rightarrow X$ to produce an arrow $A \rightarrow L(X)$. As we'll see in section \ref{sect:list:UP}, a key technical step in this process will be obtaining an intermediate arrow $g' : N \times [A_{>0}] \rightarrow X$, where $A_{>0} = (a : A \mid l_A(a)>0)$. The important property of $g'$ is that it should produce the same output as $g$, as long as it is acting on inputs $m:N$, $a:A$ such that $m < l_A(a)$.

In order to obtain the arrow $g'$ we want, we need to introduce an arrow $\IdUntil : N \times N \rightarrow N$, defined by
\begin{align*}
\IdUntil(m,n) =_{m,n} \min(m, Pn).
\end{align*}
Here, $P$ is the predecessor function, and $\min$ is defined as in \cite[Prop. 1.3]{Roman}. What's important about $\IdUntil$ is that it has the following properties.

\begin{fprop}
\label{prop:idUntil}
Let $m,n : N$ be terms in a context $C$. Then:
\begin{itemize}
\item If $0 <_C n$, then $\IdUntil(m,n) <_C n$.
\item If $m <_C n$, then $\IdUntil(m,n) =_C m$.
\end{itemize}
\end{fprop}

\begin{sproof}
For the first part, we first note that $n >_C 0$ implies $n =_C s(P(n))$ (see proposition \ref{prop:>0}). Using this and \cite[Prop. 1.3]{Roman}, we compute
\begin{align*}
s \big( \IdUntil(m,n) \big) \dotminus n
=_C s \big( \min(m,Pn) \big) \dotminus s(P(n))
=_C \min(m,Pn) \dotminus Pn.
\end{align*}
It is easy to check that $\min(x,y) \dotminus y =_{x,y} 0$ (using \cite[Prop. 1.3]{Roman}), so this gives us $\IdUntil(m,n) <_C n$, as desired.

For the second part, we note that $n >_C 0$, and so once again $n = s(P(n))$. The inequality $m <_C n$ therefore implies that $m \leq P(n)$ (since $sx \dotminus sy =_{x,y} x \dotminus y$). This fact and  \cite[Prop. 1.3]{Roman} let us compute
\begin{align*}
\IdUntil(m, n)
=_C \min(m, Pn)
=_C m \dotminus [m \dotminus Pn]
=_C m \dotminus 0
=_C m.
\end{align*}
This is all we needed to show.
\end{sproof}

We can now use the arrow $\IdUntil$ to obtain $g'$ from $g$. Essentially, $g'(m,a)$ is equal to $g(m,a)$ for $m < l_A(a)$, and for $m \geq l_A(a)$, it is just equal to $g(l_A(a) \dotminus 1, \, a)$ (which is a valid term because $l_A(a)>0$, by assumption).

\begin{fcorol}
\label{corol:idUntilEquiv}
Let $l_A : A \rightarrow N$ be an arrow, let $C = (m:N, \, a:A \mid m < l_A(a))$, and let $A_{>0} = (a : A \mid l_A(a) > 0)$. Then, for any $g : [C] \rightarrow X$, the arrow $g' : N \times [A_{>0}] \rightarrow X$ given by
\begin{align*}
g'(m,a) =_{m,a} g \big( \IdUntil(m, l_A(a)), \; a \big)
\end{align*}
is well-defined. Moreover, in the context $C$, the term $g'(m,a)$ is well-defined, and
\begin{align*}
g'(m,a) =_C g(m,a).
\end{align*}
\end{fcorol}

\begin{sproof}
We start with the first part of the statement. Let $n = \IdUntil(m, l_A(a))$; we must check that $g(n, a)$ is well-defined (remark that $a : [A_{>0}]$, but we also write $a$ for its image via the inclusion $[A_{>0}] \hookrightarrow A$). Specifically, we need to check that $(n, a)$ is a valid term of type $[C]$ by checking that $n <_{m,a} l_A(a)$. By definition of $A_{>0}$, we have $l_A(a) >_{m,a} 0$, so we can apply proposition \ref{prop:idUntil} to find that
\begin{align*}
n =_{m,a} \IdUntil(m, l_A(a))
<_{m,a} l_A(a),
\end{align*}
as desired.

For the second part, we must check that $g'(m,a)$ is well-defined. In particular, we need to check that the term $a : A$ in the context $C$ can be considered as a term of type $[A_{>0}]$. This is simple: we have $l_A(a) > m \geq 0$, which is all we need.

Finally, we must check that $g'(m,a) =_C g(m,a)$. Since $m < l_A(a)$ in $C$, we can apply proposition \ref{prop:idUntil} to get $\IdUntil(m, l_A(a)) =_C m$, and so $g'(m,a) =_C g(m,a)$, as desired.
\end{sproof}

\section{Extracting elements from lists}
\label{sect:nth}

Throughout this section, we assume we're working in a category $\mathcal{C}$ with finite limits and \param{} list objects (and therefore a \param{} natural numbers object).

As noted in section \ref{sect:outline}, the goal of this section is to construct arrows $\nth_X : E \times_N L(X) \rightarrow X$ and check that they form a natural transformation. To do so, we must first understand what we expect the arrow $\nth_X$ should be. Proposition \ref{prop:Eeq} shows that we can write
\begin{align*}
E \times_N L(X) = \{m : N, \, \ell : L(X) \mid m < \len(\ell)\},
\end{align*}
so $\nth_X$ should take a list $\ell$ and an integer $m$ with $m < \len(\ell)$, and return an element of $X$. As the name suggests, this map should output the $m^{th}$ element of $\ell$; the condition $m < \len(\ell)$ ensures that this is well-defined.

In order to construct $\nth_X$, we must first construct an arrow $\nthDef_X : X \times N \times L(X) \rightarrow X$, which acts the same way except that $\nthDef(x,m,\ell)$ returns a default element $x$ when $m \geq \len(\ell)$. This arrow can be easily defined by induction, but using it to define $\nth_X$ is more difficult, because we need to remove the default parameter. We will need to invoke hypothesis \ref{hyp:list} to resolve this issue.

\subsection{Defining nthDef}
\label{sect:nth:nthdef}

We want to define an arrow $\nthDef(x,m,\ell)$ that gives us the $m^{th}$ element of a list $\ell$ (and returns a default $x$ when $m \geq \len(\ell)$). The key idea here is that the $m^{th}$ element of $\ell$ is the zeroth element of $\tail(m,\ell)$, the list obtained by removing the leading $m$ elements of $\ell$ (see section \ref{sect:NNOlist:listars}). Therefore, we define the arrow $\nthDef_X : X \times N \times L(X) \rightarrow X$ by
\begin{align*}
\nthDef_X(x,n,\ell) =_{x,n,\ell} \zerothDef_X\Big( x, \; \tail_X\big(n, \ell\big) \Big).
\end{align*}
We often drop the subscript $X$ since this object is usually clear from context.

\subsection{Properties of nthDef}
\label{sect:nth:properties}

First, we note that $\nthDef$ is a natural transformation.

\begin{fprop}
\label{prop:nthDefNatural}
Let $f : X \rightarrow Y$ be an arrow. Then
\begin{align*}
\nthDef_Y(f(x), m, L(f)(\ell)) =_{x,m,\ell} f(\nthDef_X(x,m,\ell)).
\end{align*}
\end{fprop}

\begin{sproof}
This is a straightforward check using proposition \ref{prop:Lfunctor} and the definition of $\nthDef$.
\end{sproof}

The next two propositions tell us how $\nthDef(x,n,\ell)$ behaves depending on whether $n \geq \len(\ell)$, in which case there is no $n^{th}$ element and we should return the default $x$, or $n < \len(\ell)$, in which case the default $x$ has no impact on the outcome. Perhaps surprisingly, hypothesis \ref{hyp:list} is required for the case $n < \len(\ell)$.

\begin{fprop}
\label{prop:nthEqualsDef}
Let $x : X$, $n : N$, $\ell : L(X)$ be terms in a context $C$. If $n \geq_C \len(\ell)$, then
\begin{align*}
\nthDef(x,n,\ell) =_C x.
\end{align*}
\end{fprop}

\begin{sproof}
Since $n \geq_C \len(\ell)$ by assumption, we have $\tail(n,\ell) =_C \varnothing$ by proposition \ref{prop:lenIterTr}. So, by definition of $\nthDef$, we have
\begin{align*}
\nthDef(x,n,\ell) 
=_C \zerothDef(x, \; \tail(n,\ell))
=_C \zerothDef(x, \varnothing)
=_C x.
\end{align*}
This is what we wanted.
\end{sproof}

\begin{fprop}
\label{prop:nthNotDef}
Let $x : X$, $n : N$, $\ell : L(X)$ be terms in a context $C$, and assume \ref{hyp:list}. If $n <_C \len(\ell)$, then for any other term $x' : X$ in the context $C$,
\begin{align*}
\nthDef(x,n,\ell) =_C \nthDef(x',n,\ell).
\end{align*}
\end{fprop}

\begin{sproof}
By proposition \ref{prop:appN:calculation} from the appendix, $n<_C \len(\ell)$ implies $\len(\ell) \dotminus n >_C 0$, so
\begin{align*}
\len(\tail(n,\ell)) =_C \len(\ell) \dotminus n >_C 0
\end{align*}
by proposition \ref{prop:lenIterTr}. By the term version of hypothesis \ref{hyp:list} (proposition \ref{prop:H1term}), this inequality implies that $\tail(n,\ell) =_C y :: \ell'$ for some terms $y,\ell'$ in the context $C$. We conclude that
\begin{align*}
\nthDef(x,n,\ell) =_C \zerothDef(x, \tail(n,\ell))
=_C \zerothDef(x, y :: \ell') =_C y,
\end{align*}
and $\nthDef(x',n,\ell) =_C y$ by the same argument. This gives the desired equality.
\end{sproof}

\subsection{Removing the default parameter}
\label{sect:nth:nth}

We are now prepared to construct the arrow $\nth_X : E \times_N L(X) \rightarrow X$. Recall that, by proposition \ref{prop:Eeq}, we can use the context $C = (m : N, \, \ell : L(X) \mid m < \len(\ell))$ to talk about $E \times_N L(X)$; in this context, we would like to set $\nth_X(m,\ell) =_C \nthDef_X(m,\ell)$. However, $\nthDef_X$ always requires a default element, so to define $\nth_X$, we need a term $\deff_X : X$ in the context $C$ to serve as the default element. Then, we can set
\begin{align*}
\nth_X(m,\ell) =_C \nthDef_X (\deff_X, \, m, \, \ell).
\end{align*}
How can we obtain this $\deff_X$? This will in fact require hypothesis \ref{hyp:list}.

\begin{fprop}
\label{prop:def}
Assume \ref{hyp:list}. In the context $C = (m : N, \; \ell : L(X) \mid m < \len(\ell))$, there exists a term $\deff_X$ of type $X$.
\end{fprop}

\begin{sproof}
In the context $C$, remark that $\len(\ell) >_C m \geq_C 0$. By the term version of \ref{hyp:list} (proposition \ref{prop:H1term}), this implies that there exist terms $x : X$, $\ell' : L(X)$ in $C$ such that $\ell =_C x :: \ell'$. So, we can just take $\deff_X$ to be this term $x$. 
\end{sproof}

With this default term, our definition of $\nth_X$ is now complete. However, it's important to be aware that any further discussion of the arrow $\nth_X$ implicitly assumes the hypothesis \ref{hyp:list}.

\begin{fprop}
\label{prop:epsNT}
Assume \ref{hyp:list}. The collection of arrows $(\nth_X)_X$ is a natural transformation $\Sigma_E \circ \Delta_{\pi_2^E} \circ L_N \Rightarrow \Id_{\mathcal{C}}$.
\end{fprop}

\begin{sproof}
Let $f : X \rightarrow Y$ be an arrow of $\mathcal{C}$. We must show that the following square commutes.
\begin{equation*}
\begin{tikzcd}
E \times_N L(X) \arrow[r, "\nth_X"] \arrow[d, "\Id_E \times_N L(f)"']
& X \arrow[d, "f"] \\
E \times_N L(Y) \arrow[r, "\nth_Y"'] & Y
\end{tikzcd}
\end{equation*}
Let $C$ be the context $(m : N, \; \ell : L(X) \mid m < \len(\ell))$.
We compute:
\begin{align*}
\nth_Y (m, L(f)(\ell))
&=_C \nthDef_Y \big( \deff_Y', \; m, \; L(f)(\ell) \big)
\\
&=_C \nthDef_Y \big( f(\deff_X), \; m, \; L(f)(\ell) \big)
\\
&=_C f(\nthDef_X(\deff_X, m, \ell))
\\
&=_C f(\nth_X(m, \ell)).
\end{align*}
Here, $\deff_Y'$ is the term $\deff_Y[L(f)(\ell)/\ell]$.
For the second equality, we used proposition \ref{prop:nthNotDef}. This is allowed because we assumed \ref{hyp:list} and because $m < \len(\ell) = \len(L(f)(\ell))$ (using proposition \ref{prop:Lfunctor}). Finally, the third equality uses proposition \ref{prop:nthDefNatural}.
\end{sproof}

\pagebreak

\section{Constructing lists and computing their elements}
\label{sect:list}

Throughout this section, we assume we're working in a category $\mathcal{C}$ with finite limits and \param{} list objects (and therefore a \param{} natural numbers object).

As noted in section \ref{sect:outline}, the goal of this section is to show the existence part of the statement that $(L_N(X), \nth_X)$ is a universal morphism. That is, given arrows $l_A : A \rightarrow N$ and $g : E \times_N A \rightarrow X$, we must construct an arrow $h : A \rightarrow L(X)$ such that $l_A = \len_X \circ h$ and $g = \nth_X \circ (\Id \times_N h)$.

It is particularly difficult to construct maps into a list object because its universal property only provides inductive constructions of maps out of it. Therefore, we start this section by developing a general technique for constructing such maps. Then, we establish some properties of this technique, before finally applying it to this particular case.

\subsection{Constructing maps into a list object}

Before we start, we must introduce some notation. We write $\concat : L(X) \times L(X) \rightarrow L(X)$ for the \textit{list concatenation} arrow, defined inductively by
\begin{align*}
\varnothing \concat \ell_2 =_{\ell_2} \ell_2,
&& (x :: \ell_1) \concat \ell_2 =_{x,\ell_1,\ell_2} x :: (\ell_1 \concat \ell_2).
\end{align*}
Moreover, we'll write $[x]$ to denote the list $x :: \varnothing$.

The key idea for constructing maps into $L(X)$ is that we should use the inductive property of NNOs. Indeed, suppose we have a map $f : N \times A \rightarrow X$, which we think of as a collection of infinite sequences $(x_{n,a})_n$ indexed by $a \in A$. We can use this map to construct the map $\Seq[f] : N \times N \times A \rightarrow L(X)$, where $\Seq[f](m,n,a)$ represents $[x_{m,a}, ..., x_{m+n-1,a}]$. This is done by induction on $n$, the length of this list, as follows.
\begin{align*}
\Seq[f](m,0,a) =_{m,a} \varnothing
&& \Seq[f](m,s(n),a) =_{m,n,a} \Seq[f](m,n,a) \concat [f(m+n,a)]
\end{align*}
If we are also given a map $p : A \rightarrow N$ which represents the length of a list associated to $a$, then we can define $\List[f,p] : A \rightarrow L(X)$ by setting
\begin{align*}
\List[f,p](a) =_a \Seq[f](0, p(a), a).
\end{align*}
Intuitively, we have $\List[f,p](a) = [x_{0,a}, ..., x_{p(a)-1,a}]$. The drawback of this technique is that it requires an infinite sequence $(x_{n,a})_n$ to specify this finite list; we address this issue in section \ref{sect:list:UP}.

We end this section with some very basic facts about this construction.

\begin{fprop}
\label{prop:lenConcat}
$\len(\ell_1 \concat \ell_2) =_{\ell_1, \ell_2} \len(\ell_1) + \len(\ell_2)$.
\end{fprop}

\begin{sproof}
By induction on $\ell_1$, it suffices to show that these arrows both satisfy $h(\varnothing,\ell_2) =_{\ell_2} \len(\ell_2)$ and $h(x :: \ell_1, \; \ell_2) =_{x, \ell_1, \ell_2} 1 + h(\ell_1, \ell_2)$. This is a trivial calculation.
\end{sproof}

\begin{fprop}
\label{prop:lenList}
$\len \circ \List[f,p] = p$
\end{fprop}

\begin{sproof}
To obtain this equality, it suffices to show that that $\len(\Seq[f](m,n,a)) =_{m,n,a} n$. This is done by induction on $n$: using proposition \ref{prop:lenConcat} and noting that $\len([x]) =_x 1$, it is easy to show that these terms both satisfy the introductory equations $h(m,0,a) =_{m,a} 0$ and $h(m,sn,a) =_{m,n,a} h(m,n,a)+1$.
\end{sproof}

\begin{fprop}
\label{prop:SeqHead}
$\Seq[f](m,sn,a) =_{m,n,a} f(m,a) :: \Seq[f](sm,n,a)$.
\end{fprop}

\begin{sproof}
This is proved by induction on $n$: it is straightforward to check that both terms satisfy $h(m,0,a) =_{m,a} [f(m,a)]$ and $h(m,sn,a) =_{m,n,a} h(m,n,a) \concat [f(s(m+n),a)]$.
\end{sproof}

\subsection{The nth elements of constructed lists}

In this section, we prove an important result about the interaction between $\nth$ and $\List[f,p]$: specifically, we show with theorem \ref{thm:nthList} that the $m^{th}$ element of $\List[f,p](a)$ is $f(m,a)$, as long as $m<p(a)$.

We start with some intermediate results. The first thing we need to do is investigate the interaction between truncation and $\Seq[f]$, because $\nthDef$ is built using $\tr$, and $\List$ is built using $\Seq$. Recall that $P$ represents the predecessor function.

\begin{fprop}
\label{prop:trSeq}
Let $f : N \times A \rightarrow X$. Then $\tr(\Seq[f](m,n,a)) =_{m,n,a} \Seq[f](sm,Pn,a)$.
\end{fprop}

\begin{sproof}
By induction on $n$, it suffices to show that the arrows agree on $0$ and $s$. For the base case, it's clear that both terms reduce to $\varnothing$; for the inductive case, we use proposition \ref{prop:SeqHead}:
\begin{align*}
\tr(\Seq[f](m,sn,a)) =_{m,n,a} \tr\Big( f(m,a) :: \Seq[f](sm,n,a) \Big)
=_{m,n,a} \Seq[f](sm,n,a),
\\
\Seq[f](sm,Psn,a) =_{m,n,a} \Seq[f](sm,n,a).
\end{align*}
This is all we need to show.
\end{sproof}

From this fact, we expect that iterating the truncation map just turns into iterating the successor $s$ and predecessor $P$. However, iterating the successor just becomes addition, and iterating the predecessor just becomes subtraction. This observation gives us the following important lemma.

\begin{flem}
\label{lem:iterTrSeqList}
Let $f : N \times A \rightarrow X$. Then
\begin{align*}
\tail\Big(k, \; \Seq[f](m,n,a)\Big) =_{k,m,n,a} \Seq[f]\Big(m+k, \; n \dotminus k, \; a\Big).
\end{align*}
In particular, given $p : A \rightarrow N$, we have
\begin{align*}
\tail\Big(k, \; \List[f,p](a)\Big) =_{k,a} \Seq[f]\Big(k, \; p(a) \dotminus k, \; a\Big).
\end{align*}
\end{flem}

\begin{sproof}
We start with the first equality, which we prove by induction on $k$. To this end, we show that both arrows satisfy the equalities $h(0,m,n,a) =_{m,n,a} \Seq[f](m,n,a)$ and $h(sk,m,n,a) =_{k,m,n,a} \tr(h(k,m,n,a))$. The first equality is trivial, since $\tail(0,\ell) =_\ell \ell$, $m+0 =_m m$, and $n \dotminus 0 =_n n$. The second equality is immediate for the left-hand side function (it's the definition of $\tail$); for the right-hand side function, we use proposition \ref{prop:trSeq}.
\begin{align*}
\Seq[f]\Big(m+sk, \; n \dotminus sk, \; a \Big)
&=_{k,m,n,a} \Seq[f]\Big(s(m+k), \; P(n \dotminus k), \; a \Big)\Big)
\\
&=_{k,m,n,a} \tr \Big( \Seq[f]\Big( m+k, \; n \dotminus k, \; a\Big)\Big)
\end{align*}
This is all we needed to show for the first part. For the second part, we simply replace $m$ by $0$ and $n$ by $p(a)$, then use the definition of $\List[f,p]$.
\end{sproof}

Finally, we arrive at the desired result. Note that we need to assume hypothesis \ref{hyp:list} in order to talk about the $\nth$ arrow.

\begin{fthm}
\label{thm:nthList}
Let $f : N \times A \rightarrow X$ and $p : A \rightarrow N$, let $m:N$, $a:A$ be terms in a context $C$, and assume \ref{hyp:list}. If $m <_C p(a)$, then
\begin{align*}
\nth(m, \List[f,p](a)) =_C f(m,a).
\end{align*}
\end{fthm}

\begin{sproof}
First, we note that the term on the left is well-defined, because $\len(\List[f,p](a)) =_C p(a)$ (by proposition \ref{prop:lenList}) and $m <_C p(a)$ be assumption. Then, by the definition of $\nth$ and lemma \ref{lem:iterTrSeqList},
\begin{align*}
\nth(m, \List[f,p](a))
&=_C
\nthDef(\deff_X', m, \List[f,p](a))
\\
&=_C \zerothDef\Big(\deff_X', \; \tail\Big(m, \; \List[f,p](a)\Big)\Big)
\\
&=_C \zerothDef\Big(\deff_X', \; \Seq[f]\Big(m, p(a) \dotminus m, a \Big)\Big).
\end{align*}
(Here, $\deff_X'$ is $\deff_X$ with the appropriate substitution.)
Since $m <_C p(a)$ by assumption, we have $p(a) \dotminus m >_C 0$ by proposition \ref{prop:appN:calculation} of the appendix. Therefore, by proposition \ref{prop:>0}, $p(a) \dotminus m =_C s(k)$ for some term $k$. We substitute this into the above equality, apply proposition \ref{prop:SeqHead}, and use the definition of $\zerothDef$. We get:
\begin{align*}
\nth(m, \List[f,p](a))
&=_C \zerothDef\Big(\deff_X', \; \Seq[f](m, sk, a) \Big)
\\
&=_C \zerothDef\Big(\deff_X', \; f(m,a) :: \Seq[f](sm,k,a) \Big)
=_C f(m,a).
\end{align*}
This is what we wanted to show.
\end{sproof}

\subsection{The universal property arrow}
\label{sect:list:UP}

The construction $\List[f,p]$ allows us to construct an arrow $A \rightarrow L(X)$ based on two arrows $f : N \times A \rightarrow X$ and $p : A \rightarrow N$. However, the goal stated at the beginning of this chapter is to construct $h : A \rightarrow L(X)$ from $g : E \times_N A \rightarrow X$ and $l_A : A \rightarrow N$, so we need to adjust this technique.

Recall from proposition \ref{prop:Eeq} that $E \times_N A = \{m : N, \, a : A \mid m < l_A(a)\}$. This means that the arrow $g : E \times_N A \rightarrow X$ only gives us finite sequences $(x_{0,a}, ..., x_{l_A(a)-1,a})$ instead of the infinite sequences used in $\List[f,p]$. We account for this discrepancy by using two tricks.
\begin{itemize}
\item First, if $l_A(a)>0$, then we can extend the finite sequence to an infinite one simply by repeating the last element; we will use corollary \ref{corol:idUntilEquiv} to do this formally.

\item Second, we note that the first trick doesn't work if $l_A(a)=0$, because then $g$ doesn't give us any elements of $X$. However, $l_A(a)=0$ just means that we should set $h(a)=\varnothing$. So, we have to use hypothesis \ref{hyp:NNO} to split into two cases: $l_A(a)=0$ and $l_A(a)>0$.
\end{itemize}
We employ these techniques in the following proof. Note that hypothesis \ref{hyp:list} is required to talk about $\nth_X$ (see section \ref{sect:nth:nth}).

\begin{fthm}
\label{thm:UPexistence}
Assume \ref{hyp:list} and \ref{hyp:NNO}. Given objects $X$, $A$ and arrows $l_A : A \rightarrow N$ and $g : E \times_N A \rightarrow X$, there exists an arrow $h : A \rightarrow L(X)$ such that $l_A = \len_X \circ h$ and $g = \nth_X \circ (\Id_E \times_N h)$.
\end{fthm}

\begin{sproof}
Hypothesis \ref{hyp:NNO} tells us (via proposition \ref{prop:coproduct}) that we have the coproduct $A = [A_0] + [A_{>0}]$, where $A_0 = (a : A \mid l_A(a)=0)$ and $A_{>0} = (a : A \mid l_A(a) > 0)$. This lets us define $h : A \rightarrow L(X)$ ``by cases" on each part of this coproduct.

In the case $l_A(a)=0$, we simply set $h(a) =_{A_0} \varnothing$. For the case $l_A(a) > 0$, note that the domain of $g$ is $[C]$, where $C = (m : N, \, a : A \mid m < l_A(a))$. Then corollary \ref{corol:idUntilEquiv} tells us that there is an arrow $g' : N \times [A_{>0}] \rightarrow X$ such that $g'(m,a) =_C g(m,a)$, and we can define $h$ in the $A_{>0}$ case to be
\begin{align*}
h(a) =_{A_{>0}} \List[g', l_A'](a),
\end{align*}
where $l_A' : [A_{>0}] \rightarrow N$ is just given by $l_A'(a) =_{A_{>0}} l_A(a)$. This completes the definition of $h$; we must now check that it satisfies the required equalities.

First, we check that $\len(h(a)) =_a l_A(a)$. This can be done by checking equality on both parts of the coproduct. In the $A_0$ case, we have $\len(h(a)) =_{A_0} \len(\varnothing) =_{A_0} 0 =_{A_0} l_A(a)$, and in the $A_{>0}$ case, we have
\begin{align*}
\len(h(a)) =_{A_{>0}} \len(\List[g', l_A'](a)) =_{A_{>0}} l_A'(a) =_{A_{>0}} l_A(a)
\end{align*}
by using proposition \ref{prop:lenList}.

Second, we check that $\nth_X(m, h(a)) =_C g(m,a)$. We could check on both parts of the coproduct, but this would be redundant: the condition $m<l_A(a)$ in $C$ implies that we are already in the $A_{>0}$ case. Formally, performing the substitution $a \mapsto a$ in the equation $h(a) =_{A_{>0}} \List[g',l_A'](a)$ gives us the equality $h(a) =_C \List[g', l_A'](a)$, and this substitution is valid because $l_A(a) > 0$ in $C$. Using this equality and theorem \ref{thm:nthList}, we compute
\begin{align*}
\nth_X(m,h(a)) 
=_C \nth_X(m, \List[g', l_A'](a))
=_C g'(m, a)
=_C g(m,a),
\end{align*}
which is what we wanted.
\end{sproof}

\section{Equality of lists}
\label{sect:equality}

Throughout this section, we assume we're working in a category $\mathcal{C}$ with finite limits and \param{} list objects (and therefore a \param{} natural numbers object).

As noted in section \ref{sect:outline}, the goal of this section is to show the uniqueness part of the statement that $(L_N(X), \nth_X)$ is a universal morphism. That is, if we have two arrows $h_1, h_2 : A \rightarrow L(X)$ such that $l_A = \len_X \circ h_i$ and $g = \nth_X \circ (\Id \times_N h_i)$ for $i=1,2$, then $h_1 = h_2$ (here, $l_A : A \rightarrow N$ and $g : E \times_N A \rightarrow X$). Intuitively, what we need to show is that if two lists have the same length and the same elements, then they are equal.

To simplify this problem, we can start by using the arrow $\nthDef$ instead of $\nth$. So, given two arrows $h_1, h_2 : A \rightarrow L(X)$, we'd like to show the following statement:
\begin{center}
If $\len(h_1(a)) =_a \len(h_2(a))$ and $\nthDef(x,n,h_1(a)) =_{x,n,a} \nthDef(x,n,h_2(a))$,\\then $h_1(a) =_a h_2(a)$.
\end{center}
It turns out that hypotheses \ref{hyp:list} and \ref{hyp:NNO} are necessary to prove this statement. Without them, the best we can do is show that $h_1(a) =_{x,a} h_2(a)$, where $x:X$; this is the content of theorem \ref{thm:listEqNoExt}.

This section will therefore be divided into two subsections. In the first, we prove theorem \ref{thm:listEqNoExt}, which doesn't require \ref{hyp:list} or \ref{hyp:NNO}. In the second, we use these hypotheses to improve this result: first, we change the conclusion to be $h_1 = h_2$, and then we adjust our hypotheses so they involve $\nth$ instead of $\nthDef$.

\subsection{Equality without extra hypotheses}

Our goal for this section is to prove theorem \ref{thm:listEqNoExt}, below. Essentially, we want to prove that two lists $h_1, h_2 : A \rightarrow L(X)$ are equal if they have the same length and the same elements. The idea of the proof is to go by induction on their length, and show that they are ``built up in the same way".

More specifically, the strategy is as follows. Given arrows $h_1, h_2 : A \rightarrow L(X)$, we form arrows $H_1, H_2 : N \times A \rightarrow L(X)$ where $H_i(k,a)$ is the list $h_i(a)$ with all but the last $k$ elements removed. We then want to show that $H_1 = H_2$; if this is true, then by setting $k = \len(h_i(a))$, we get $h_1(a) = h_2(a)$, as desired. To show that $H_1 = H_2$, we go by induction on $k$: the base case is clear since $H_1(0,a) = H_2(0,a) = \varnothing$, but the inductive step is more tricky.

For the inductive step, we must analyse the transition from having the last $k$ elements of $h_i$ to having the last $k+1$ elements of $h_i$. What we're doing is ``adding back" an element to $h_i$; specifically, the one in position $\len(h_i) \dotminus k$. So, we'd expect to have
\begin{align*}
H_i(k+1,a) = \nthDef\big(\len(h_i(a)) \dotminus k, h_i(a)\big) :: H_i(k,a).
\end{align*}
If we're assuming that $\nthDef$ and $\len$ agree for $h_1$ and $h_2$, then this gives us a recurrence relation that both $H_1$ and $H_2$ satisfy, and we're done. However, there's a problem: $\nthDef$ must include a default parameter. For this reason, $H_i$ must also include a default parameter, making it of the form $H_i(x,k,a)$. This is why we can only conclude $h_1(a) =_{x,a} h_2(a)$.

Formalizing this proof idea requires several steps, and we will make heavy use of the if-then-else notation from section \ref{sect:NNOlist:ITE}. First, we must have a better understanding of the arrows $H_i$. These arrows will involve the mapping $(x,k,\ell) \mapsto \tail(\len(\ell) \dotminus k, \ell)$, which we want to show interacts nicely with appending an $n^{th}$ element of $\ell$. These first two propositions give us more information about this idea.

\begin{fprop}
We have
\begin{align*}
\tail(m,\ell) =_{x,m,\ell} \left\{\begin{matrix}
\varnothing & \text{if } \len(\ell) \leq m \\
\nthDef(x,m,\ell) :: \tail(sm,\ell) & \text{else}
\end{matrix}\right. .
\end{align*}
\end{fprop}

\begin{sproof}
We first prove by induction on $\ell$ that
\begin{align*}
\ell =_{x,\ell} \left\{\begin{matrix}
\varnothing & \text{if } \len(\ell) = 0 \\
\zerothDef(x, \ell) :: \tr(\ell) & \text{else}
\end{matrix}\right. .
\end{align*}
If $\ell = \varnothing$, then $\len(\ell) = 0$, so the right side reduces to $\varnothing$. If $\ell = y :: \ell'$, then $\len(\ell) = s(\len(\ell'))$, and the right side reduces to $\zerothDef(x, y :: \ell') :: \tr(y :: \ell') = y :: \ell' = \ell$.

Now, in the above equality, we replace $\ell$ by $\tail(m,\ell)$. The else case simplifies (using the definition of $\nthDef$) to
\begin{align*}
\zerothDef(x, \tail(m,\ell)) :: \tr(\tail(m,\ell)) = \nthDef(x,m,\ell) :: \tail(sm,\ell).
\end{align*}
The term $\len(\tail(m,\ell))$ in the ``if" condition simplifies (by proposition \ref{prop:lenIterTr}) to $\len(\ell) \dotminus m$, so the condition is equivalent to $\len(\ell) \leq m$. Thus, we get the desired equality.
\end{sproof}

\begin{fprop}
We have
\begin{align*}
\tail(Pm,\ell) =_{x,m,\ell} \left\{\begin{matrix}
\tail(m,\ell) & \text{if } m=0 \\
\varnothing & \text{else if } \len(\ell) < m \\
\nthDef(x,Pm,\ell) :: \tail(m,\ell) & \text{else}
\end{matrix}\right. .
\end{align*}
\end{fprop}

\begin{sproof}
We prove this by induction on $m$. If $m=0$, then the left side is just $\tail(0,\ell)$. On the right side, the condition $m=0$ is true, so it reduces to $\tail(0,\ell)$, as desired.

On the other hand, suppose $m=sn$. On the right hand side, the outer condition becomes $sn=0$, which is false, so it reduces to the else-if. So, the equality we want to prove becomes
\begin{align*}
\tail(Psn,\ell) =_{x,n,\ell} \left\{\begin{matrix}
\varnothing & \text{if } \len(\ell) < sn \\
\nthDef(x,Psn,\ell) :: \tail(sn,\ell) & \text{else}
\end{matrix}\right. .
\end{align*}
Clearly $Psn=n$, and the condition $\len(\ell) < sn$ is by definition equivalent to $s(\len(\ell)) \dotminus sn = 0$. But $sa \dotminus sb =_{a,b} a \dotminus b$, so this condition is just $\len(\ell) \dotminus n = 0$, i.e. $\len(\ell) \leq n$. Thus, the equality we want to prove is precisely what we have proved in the previous proposition.
\end{sproof}

In the next lemma, we show that an arrow $H$ (which will act like the arrows $H_1, H_2$ we mentioned) satisfies a recurrence relation like the one we described in the proof outline.

\begin{flem}
Consider the arrow $A : X \times N \times L(X) \times L(X) \rightarrow L(X)$ given by
\begin{align*}
A(x,k,\ell,L) =_{x,k,\ell,L} \left\{
\begin{matrix}
L & \text{if } \len(\ell) \leq k \\
\nthDef\big(x, P(\len(\ell)\dotminus k), \ell\big) :: L
& \text{else}
\end{matrix}
\right. .
\end{align*}
Next, consider $H : X \times N \times L(X) \rightarrow L(X)$ given by
\begin{align*}
H(x,k,\ell) =_{x,k,\ell} \tail(\len(\ell)\dotminus k, \ell).
\end{align*}
Intuitively, $H$ takes $\ell$ and removes all but the last $k$ elements. Then,
\begin{align*}
H(x,s(k),\ell) =_{x,k,\ell} A(x,k,\ell,H(x,k,\ell)).
\end{align*}
\end{flem}

\begin{sproof}
Using the previous proposition, we compute
\begin{align*}
& H(x,sk,\ell) 
\\
&=_{x,k,\ell} \tail(\len(\ell) \dotminus sk, \ell)
\\
&=_{x,k,\ell} \tail(P(\len(\ell) \dotminus k), \ell)
\\
&=_{x,k,\ell} 
\left\{\begin{matrix}
\tail(\len(\ell) \dotminus k,\ell) & \text{if } \len(\ell) \dotminus k=0 \\
\varnothing & \text{else if } \len(\ell) < \len(\ell) \dotminus k \\
\nthDef(x,P(\len(\ell) \dotminus k),\ell) :: \tail(\len(\ell) \dotminus k,\ell)  & \text{else}
\end{matrix}\right. .
\end{align*}
Note that the top condition, $\len(\ell) \dotminus k = 0$, is the same as $\len(\ell) \leq k$. Moreover, the else-if condition is just false, since $s(\len(\ell)) \dotminus (\len(\ell) \dotminus k)$ is greater than zero by proposition \ref{prop:appN:natCalc1} of the appendix. Thus, the else-if-else reduces to just the else case. Finally, we can just replace the term $\tail(\len(\ell) \dotminus k, \ell)$ by $H(x,k,\ell)$. With all this consideration, the above expression reduces to
\begin{align*}
H(x,sk,\ell) = \left\{\begin{matrix}
H(x,k,\ell) & \text{if } \len(\ell) \leq k \\
\nthDef(x, P(\len(\ell) \dotminus k), \ell) :: H(x,k,\ell) & \text{else}
\end{matrix}\right. ,
\end{align*}
which is just $A(x,k,\ell, H(x,k,\ell))$, as desired.
\end{sproof}

Finally, we reach the desired result. As noted previously, this theorem does not allow us to conclude $h_1(a) =_a h_2(a)$, i.e. $h_1 = h_2$. It says that $h_1 \circ \pi_A = h_2 \circ \pi_A$ as maps $X \times A \rightarrow L(X)$; we can't yet get rid of the parameter $X$.

\begin{fthm}
\label{thm:listEqNoExt}
Let $h_1, h_2 : A \rightarrow L(X)$. Suppose that 
\begin{gather*}
\len(h_1(a)) =_a \len(h_2(a)),
\\
\nthDef(x, m, h_1(a)) =_{x,m,a} \nthDef(x, m, h_2(a)).
\end{gather*}
Then $h_1(a) =_{x,a} h_2(a)$, where $x$ is a variable of type $X$.
\end{fthm}

\begin{sproof}
For $i=1,2$, we define the map $A_i : X \times N \times A \times L(X) \rightarrow L(X)$ by setting $A_i(x,k,a,L) =_{x,k,a,L} A(x,k,h_i(a),L)$, where $A$ is as in the previous lemma. For $i=1,2$, define $H_i : X \times N \times A \rightarrow L(X)$ by $H_i(x,k,a) =_{x,k,a} H(x,k,h_i(a))$, where $H$ is as in the previous lemma. By the previous lemma, we know that $H_i(x,sk,a) =_{x,k,a} A_i(x,k,a,H_i(x,k,a))$.

We claim that $H_1 = H_2$; we'll prove this by induction on $k$. In the base case, we have $H_1(x,0,a) =_{x,a} \varnothing =_{x,a} H_2(x,0,a)$; indeed, for $i=1,2$, we compute
\begin{align*}
H_i(x,0,a) =_{x,a} H(x,0,h_i(a)) &=_{x,a} \tail(\len(h_i(a)) \dotminus 0, h_i(a))
\\
&=_{x,a} \tail(\len(h_i(a)), h_i(a))
=_{x,a} \varnothing.
\end{align*}
The last equality follows from combining propositions \ref{prop:lenIterTr} and \ref{prop:listLenDecomp}. (The length of this list is $\len(h_i(a)) \dotminus \len(h_i(a)) =_{x,a} 0$, so the list is empty.)

For the inductive step, we note that $A_1 = A_2$. Indeed, we have
\begin{align*}
A_i(x,k,a,L) =_{x,k,a,L} \left\{\begin{matrix}
L & \text{if } \len(h_i(a)) \leq k \\
\nthDef(x, P(\len(h_i(a)) \dotminus k), h_i(a)) :: L
& \text{else}
\end{matrix}\right. ,
\end{align*}
so the equalities assumed for this theorem give us $A_1 = A_2$. Writing $A_\times$ for $A_1 = A_2$, we therefore have $H_i(x,sk,a) = A_\times(x,k,a,H_i(x,a,k))$ for $i=1,2$.

So, $H_1, H_2$ satisfy the same defining equations, and we conclude by induction that $H_1 = H_2$. Now, note that
\begin{align*}
H_i(x, \len(h_i(a)), a) =_{x,a} \tail(0, h_i(a)) =_{x,a} h_i(a),
\end{align*}
so $H_1 = H_2$ and $\len(h_1(a)) =_a \len(h_2(a))$ let us conclude $h_1(a) =_{x,a} h_2(a)$, as desired.
\end{sproof}

\subsection{Equality with extensivity hypotheses}

In this section, we use hypotheses \ref{hyp:list} and \ref{hyp:NNO} to supplement theorem \ref{thm:listEqNoExt}. First, we show that the conclusion can be changed from $h_1(a) =_{x,a} h_2(a)$ to $h_1 = h_2$.

\begin{flem}
\label{lem:listEqExtLast}
Assume \ref{hyp:list} and \ref{hyp:NNO}. If $h_1, h_2 : A \rightarrow L(X)$ are arrows which satisfy the equations $\len(h_1(a)) =_a \len(h_2(a))$ and $h_1(a) =_{x,a} h_2(a)$ ($x$ is a variable of type $X$), then $h_1 = h_2$.
\end{flem}

\begin{sproof}
Let $l_A = \len_X \circ h_1 = \len_X \circ h_2$. Let $D_0$ be the context $(a:A \mid l_A(a)=0)$ and let $D_{>0}$ be the context $(a:A \mid l_A(a) > 0)$.
By proposition \ref{prop:coproduct} (which uses \ref{hyp:NNO}), we have $A = [D_0] + [D_{>0}]$. Therefore, to show $h_1 = h_2$, it suffices to show that
\begin{align*}
h_1(a) =_{D_0} h_2(a)
&& \text{and} &&
h_1(a) =_{D_{>0}} h_2(a).
\end{align*}

For the first equality, we claim that both terms equal $\varnothing$. Indeed, in the context $D_0$, we have $\len(h_1(a)) =_{D_0} l_A(a) =_{D_0} 0$, and this implies that $h_1(a) =_{D_0} \varnothing$ by proposition \ref{prop:listLenDecomp}. The same goes for $h_2$.

For the second equality, we need the existence of a term $\deff' : X$ in the context $D_{>0}$. Indeed, if such a term exists, then by substituting $x \mapsto \deff'$, $a \mapsto a$ in the equation $h_1(a) =_{x,a} h_2(a)$, we get the desired equality
\begin{align*}
h_1(a) =_{D_{>0}} h_2(a).
\end{align*}
To obtain $\deff'$, we use proposition \ref{prop:def} (which requires \ref{hyp:list}) to get a term $\deff : X$ in the context $(m:N, \, \ell : L(X) \mid m < \len(\ell))$. Then, since $0 <_{D_{>0}} l_A(a) =_{D_{>0}} \len(h_1(a))$, we can perform a substitution $m \mapsto 0$, $\ell \mapsto h_1(a)$ to get the term $\deff' = \deff[0/m, \, h_1(a)/\ell]$ of type $X$ in the context $D_{>0}$, as desired.
\end{sproof}

Next, we show that the hypotheses of theorem \ref{thm:listEqNoExt} can be changed so that we use $\nth$ instead of $\nthDef$ (recall that assumption \ref{hyp:list} is necessary to talk about $\nth_X$).

\begin{flem}
\label{lem:listEqExtFst}
Assume \ref{hyp:list} and \ref{hyp:NNO}. Let $h_1, h_2 : A \rightarrow L(X)$. Suppose that
\begin{align*}
\len \circ h_1 &= \len \circ h_2,
\end{align*}
and denote this arrow $l_A$. Then, in the context $C = (m : N, \; a : A \mid m < l_A(a))$, suppose that
\begin{align*}
\nth(m, h_1(a)) =_C  \nth(m, h_2(a)).
\end{align*}
Then $\nthDef(x, m, h_1(a)) =_{x,m,a} \nthDef(x, m, h_2(a))$.
\end{flem}

\begin{sproof}
Let $D = (x:X, m:N, a:A)$, and let $D_0 = (D \mid m < l_A(a))$, $D_1 = (D \mid m \geq l_A(a))$. Since we're assuming \ref{hyp:NNO}, we can use corollary \ref{corol:coproductLeq} to conclude that $[D] = [D_0] + [D_1]$. Therefore, to show that $\nthDef(x, m, h_1(a)) =_D \nthDef(x, m, h_2(a))$, it suffices to show the following.
\begin{align*}
\nthDef(x, m, h_1(a)) =_{D_0} \nthDef(x, m, h_2(a))
\\
\nthDef(x, m, h_1(a)) =_{D_1} \nthDef(x, m, h_2(a))
\end{align*}
We start with the second equality, since it is the easiest. Indeed, since we have $m \geq_{D_1} l_A(a) =_{D_1} \len(h_1(a))$, proposition \ref{prop:nthEqualsDef} tells us that
\begin{align*}
\nthDef(x, m, h_1(a)) =_{D_1} x.
\end{align*}
The same is true if we replace $h_1$ by $h_2$, so the equality follows.

For the first equality, we note that $\nth(m, h_1(a)) =_{D_0} \nth(m, h_2(a))$ by assumption (we can ``add $x:X$ to the context" by performing a substitution). To establish the desired equality, then, it suffices to show that
\begin{align*}
\nthDef(x, m, h_1(a)) =_{D_0} \nth(m, h_1(a)),
\end{align*}
and the same for $h_2$ (which is done the same way). We finish the proof by establishing this equality.
\begin{align*}
\nth(m, h_1(a))
&=_{D_0} \nthDef(\deff', m, h_1(a))
\\&=_{D_0} \nthDef( x, m, h_1(a))
\end{align*}
There are two things to note here. First, the second equality is by proposition \ref{prop:nthNotDef} (which uses \ref{hyp:list}), using the hypothesis $m < l_A(a) = \len(h_1(a))$ from the context $D_0$. Second, $\deff'$ is the term $\deff$ used to define $\nth$, but with the substitution $\ell \mapsto h_1(a)$ (it is originally in the context $(m:N, \ell:L(X) \mid m < \len(\ell))$).
\end{sproof}

Putting this all together, we get the following result, which is the uniqueness part of showing that $(L_N(X), \nth_X)$ is a universal morphism.

\begin{fcorol}
\label{corol:listEqExt}
Assume \ref{hyp:list} and \ref{hyp:NNO}. If $h_1, h_2 : A \rightarrow L(X)$ are such that $\len_X \circ h_1 = \len_X \circ h_2$ and $\nth_X \circ (\Id_E \times_N h_1) = \nth_X \circ (\Id_E \times_N h_2)$, then $h_1 = h_2$.
\end{fcorol}

\begin{sproof}
Apply lemma \ref{lem:listEqExtFst}, then theorem \ref{thm:listEqNoExt}, then lemma \ref{lem:listEqExtLast}. Remark that, if we set $l_A = \len_X \circ h_1 = \len_X \circ h_2$, then in the context $C = (m : N, a : A \mid m < l_A(a))$, we have $\big( \nth_X \circ (\Id_E \times_N h) \big)(m, \ell) =_C \nth_X(n, h(a))$ by fact \ref{fact:Idxf}.
\end{sproof}

\pagebreak
\bibliographystyle{wmaainf}
\bibliography{sources}

\pagebreak
\begin{appendices}
\section{Arithmetic for natural numbers objects}
\label{app:NNO}

In this appendix, we develop various facts about the arithmetic and order relations on a natural numbers object. Recall that the operations $+, \bullet, P, \dotminus$ are defined by induction with the following equations.
\begin{align*}
x+0 &=_x x &
x\bullet 0 &=_x 0 &
P(0) &=_{\varnothing} 0 &
x \dotminus 0 &=_x x
\\
x+sy &=_{x,y} s(x+y) &
x \bullet sy &=_{x,y} (x \bullet y) + x &
P(sy) &=_y y &
x \dotminus sy &=_{x,y} P(x \dotminus y)
\end{align*}
Moreover, given terms $m,n : N$ in a context $C$, we write $m \leq_C n$ for $m \dotminus n =_C 0$, and we write $m <_C n$ for $s(m) \dotminus n =_C 0$.

We will borrow the function $|x,y| =_{x,y} (x \dotminus y) + (y \dotminus x)$ from Roman's paper on natural numbers objects \cite{Roman} (originally from \cite{Goodstein}), which has the property that for $f,g : X \rightarrow N$, $f=g$ if and only if $|f,g| = 0_X$ (corollary 1.4 in his paper). Equivalently, given any two terms $x,y : N$ in a context $C$, we have $x =_C y$ if and only if $|x,y| =_C 0$. We will also use Roman's $\min$ and $\max$ functions, defined as $\min(x,y) =_{x,y} x \dotminus (x \dotminus y) =_{x,y} y \dotminus (y \dotminus x)$ and $\max(x,y) =_{x,y} x + (y \dotminus x) =_{x,y} y + (x \dotminus y)$ (see his proposition 1.3).

\subsection{Basic arithmetic}

\begin{fprop}
\label{prop:appN:basic}
We have the following equalities.
\begin{enumerate}
\item $a \dotminus (b+c) =_{a,b,c} (a \dotminus b) \dotminus c$
\item $|\max(a,b), b| =_{a,b} a \dotminus b$
\item $(a \dotminus b) \bullet (b \dotminus a) =_{a,b} 0$
\end{enumerate}
\end{fprop}

\begin{sproof}
For part 1, we go by induction on $c$. In the base case $c=0$, both sides clearly equal $a \dotminus b$. For the induction step, we compute
\begin{align*}
a \dotminus (b+s(c))
=_{a,b,c} a \dotminus s(b+c)
=_{a,b,c} P(a \dotminus (b+c)),
\\
(a \dotminus b) \dotminus s(c)
=_{a,b,c} P((a \dotminus b) \dotminus c).
\end{align*}
So, the equality follows by induction.

For part 2, since $|\max(a,b),b| =_{a,b} (\max(a,b) \dotminus b) + (b \dotminus \max(a,b))$, it suffices to show $\max(a,b) \dotminus b =_{a,b} a \dotminus b$ and $b \dotminus \max(a,b) =_{a,b} 0$. The latter equality is an easy computation with part 1:
\begin{align*}
b \dotminus \max(a,b)
=_{a,b} b \dotminus (a + (b \dotminus a))
=_{a,b} (b \dotminus a) \dotminus (b \dotminus a)
=_{a,b} 0.
\end{align*}
For the former equality, we use corollary 2.2(e) from Roman's paper \cite{Roman}:
\begin{align*}
\max(a,b) \dotminus b
=_{a,b} (a + (b \dotminus a)) \dotminus b
&=_{a,b} (a \dotminus b) + ((b \dotminus a) \dotminus (b \dotminus a))
\\&=_{a,b} (a \dotminus b) + 0
=_{a,b} a \dotminus b.
\end{align*}

For part 3, we will use the induction principle provided by proposition 3.1(c) in Roman's paper \cite{Roman}. This principle states that, given $f : N^2 \rightarrow N$, if $f(a+sb,b) =_{a,b} 0$, then $(a \dotminus b) \bullet f(a,b) =_{a,b} 0$. This conclusion is what we want with $f(a,b) =_{a,b} b \dotminus a$, so we compute
\begin{align*}
f(a+sb,b) =_{a,b} b \dotminus (a+sb)
=_{a,b} (b \dotminus sb) \dotminus a
=_{a,b} 0 \dotminus a =_{a,b} 0.
\end{align*}
This is all we needed to show.
\end{sproof}

\subsection{Order facts}

\begin{fprop}
\label{prop:appN:leqEquiv}
Let $a,b : N$ be terms in a context $C$. Then TFAE:
\begin{enumerate}
\item $a \leq_C b$
\item $\max(a,b) =_C b$
\item $\min(a,b) =_C a$
\item There exists a term $x : N$ in the context $C$ such that $a+x =_C b$
\item There exists a term $x : N$ in the context $C$ such that $a =_C b \dotminus x$.
\end{enumerate}
\end{fprop}

\begin{sproof}
For this proof, we show that $1 \Rightarrow 2 \Rightarrow 4 \Rightarrow 1$ and $1 \Rightarrow 3 \Rightarrow 5 \Rightarrow 1$. All the equalities here will take place in the context $C$, so we will omit the subscripts.

($1 \Rightarrow 2$). If $a \leq b$, then $a \dotminus b = 0$. Using proposition \ref{prop:appN:basic}(2), we find that $|\max(a,b),b| = a \dotminus b = 0$, which implies $\max(a,b) = b$, as desired.

($2 \Rightarrow 4$). Suppose $\max(a,b) = b$, and set $x = b \dotminus a$. Then $a+x = \max(a,b) = b$.

($4 \Rightarrow 1$). Assume $a+x=b$. Then, using \ref{prop:appN:basic}(1), $a \dotminus b = a \dotminus (a+x) = (a \dotminus a) \dotminus x = 0 \dotminus x = 0$, which is equivalent to $a \leq b$.

($1 \Rightarrow 3$). Suppose $a \leq b$, i.e. $a \dotminus b = 0$. By definition, $\min(a,b) = a \dotminus (a \dotminus b) = a \dotminus 0 = a$.

($3 \Rightarrow 5$). Suppose $\min(a,b)=a$, and set $x = b \dotminus a$. Then $b \dotminus x = \min(b,a) = \min(a,b) = a$.

($5 \Rightarrow 1$). Assume $a = b \dotminus x$. By Roman's proposition 1.3(f), $a \dotminus b = (b \dotminus x) \dotminus b = (b \dotminus b) \dotminus x = 0 \dotminus x = 0$, which is equivalent to $a \leq b$.
\end{sproof}

\begin{fprop}
\label{prop:appN:orderStandard}
For any context $C$, the relation $\leq_C$ is a partial order on terms of type $N$ in the context $C$. That is:
\begin{enumerate}
\item $x \leq_C x$;
\item If $x \leq_C y$ and $y \leq_C x$, then $x =_C y$;
\item If $x \leq_C y$ and $y \leq_C z$, then $x \leq_C z$.
\end{enumerate}
Moreover, if $x \leq_C y$, then $x+z \leq_C y+z$, $x \bullet z \leq_C y \bullet z$, and $x \dotminus z \leq_C y \dotminus z$.
\end{fprop}

\begin{sproof}
We will omit the $C$ subscript for clarity.
For 1, see Roman, proposition 1.3(b). For 2, note that $x \leq y$ and $y \leq x$ imply $|x,y| = (x \dotminus y) + (y \dotminus x) = 0 + 0 = 0$. By Roman, corollary 1.4, we get $x=y$.

For 3, we use proposition \ref{prop:appN:leqEquiv} to conclude that there exist $a,b : N$ such that $x+a = y$ and $y+b = z$. Then $x+(a+b) = (x+a)+b = y+b = z$, so $x \leq z$, as desired.

For the last part, assume that $x \leq y$, so $y = x+k$ for some $k$ (by proposition \ref{prop:appN:leqEquiv}). For the addition inequality, $y+z = (x+k)+z = (x+z)+k$, so $x+z \leq y+z$. For the multiplication, $y \bullet z = (x+k) \bullet z = x \bullet z + k \bullet z$, so $x \bullet z \leq y \bullet z$. For the subtraction, we apply Roman's corollary 2.2(e) to find
\begin{align*}
y \dotminus z
= (x+k) \dotminus z
= (x \dotminus z) + (k \dotminus (z \dotminus x)),
\end{align*}
so $x \dotminus z \leq y \dotminus z$.
\end{sproof}

\subsection{Facts needed specifically for this paper}

\begin{fprop}
\label{prop:appN:calculation}
Let $a,b : N$ be terms in a context $C$. Then $a <_C b$ if and only if $b \dotminus a >_C 0$.
\end{fprop}

\begin{sproof}
We omit the $C$ subscript for brevity. Remark that the proposition is the same as saying $b \geq a+1$ if and only if $b \dotminus a \geq 1$.

First, suppose $b \geq a+1$. Then $b \dotminus a \geq (a+1) \dotminus a = 1$ by proposition \ref{prop:appN:orderStandard} and Roman's proposition 1.3(c).

On the other hand, suppose $b \dotminus a \geq 1$. Then, using propositions \ref{prop:appN:basic}(3) and \ref{prop:appN:orderStandard},
\begin{align*}
a \dotminus b = (a \dotminus b) \bullet 1 \leq (a \dotminus b) \bullet (b \dotminus a) = 0,
\end{align*}
so $a \dotminus b = 0$. Thus,
\begin{align*}
b = b + (a \dotminus b) = a + (b \dotminus a) \geq a+1.
\end{align*}
(The second equality is the definition of $\max$; see Roman's proposition 1.3(e).)
\end{sproof}

\begin{fprop}
\label{prop:appN:natCalc1}
We have $sx \dotminus (x \dotminus y) =_{x,y} s(x \dotminus (x \dotminus y))$.
\end{fprop}

\begin{sproof}
By Roman's corollary 2.2(e),
\begin{align*}
sx \dotminus (x \dotminus y)
&= (x+1) \dotminus (x \dotminus y)
\\&= (x \dotminus (x \dotminus y)) + (1 \dotminus ((x \dotminus y) \dotminus x))
\\&= (x \dotminus (x \dotminus y)) + (1 \dotminus 0)
\\&= s(x \dotminus (x \dotminus y)).
\end{align*}
Note that $(x \dotminus y) \dotminus x = (x \dotminus x) \dotminus y = 0 \dotminus y = 0$ by Roman's proposition 1.3(e).
\end{sproof}

\end{appendices}

\end{document}